\documentclass[12pt]{article}

\usepackage{latexsym}
\usepackage{amsmath}
\usepackage{amssymb}
\usepackage{amsfonts}
\usepackage{color}

\usepackage{chngcntr}
\numberwithin{equation}{section}

\usepackage{setspace}

\begin{document}
\parskip 0.2 cm

\newtheorem{Theorem}{Theorem}[section]
\newtheorem{Proposition}[Theorem]{Proposition}
\newtheorem{Remark}[Theorem]{Remark}
\newtheorem{Lemma}[Theorem]{Lemma}
\newtheorem{Corollary}[Theorem]{Corollary}		

\def \o{\over}
\def \s{\sigma}
\def \G{\Gamma}
\def \g{\gamma}
\def \d{\delta}
\def\gm{\gamma(dy,du)}
\def\xm{\xi(dy,du)}
\def \a{\alpha}
\def \l{\lambda}
\def \disp{\displaystyle}
\def \up{\uparrow}
\def \U{{\cal U}}
\def \K{{\cal K}}
\def \P{{\cal P}}
\def \D{{\cal D}}
\def \M{{\cal M}}
\def \Mp{{\cal M_+}}
\def \ve{\varepsilon}
\def\reals{I\!\!R}
\def \O{\Omega}
\def \ph{\varphi}
\def\by{\bar y}
\def\bu{\bar u}
\def\bp{\bar \psi}
\def\be{\bar \eta}
\def\hf{\hfill{$\Box$}}
\def\T{\Theta}

\def\A{\mathcal{A}}
\def\M{\mathcal{M}}
\def\B{\mathcal{B}}

\title{LP Formulations of  Discrete Time Long-Run Average Optimal Control Problems: The Non-Ergodic Case}

\author{Vivek S. Borkar\thanks{Department of Electrical Engineering, Indian Institute of Technology Bombay,
              Powai, Mumbai 400076, India, borkar.vs@gmail.com; the work of this author was supported by a J.\ C.\ Bose Fellowship from the Government of India},    Vladimir Gaitsgory\thanks{Department of Mathematics and Statistics, Macquarie University,
              Sydney, NSW 2109, Australia,
              vladimir.gaitsgory@mq.edu.au; the work of this author was supported by the Australian Research Council Discovery Grants DP130104432} and Ilya Shvartsman\thanks{Department of Mathematics and Computer Science, Penn State Harrisburg, Middletown, PA 17057, USA, ius13@psu.edu} }

\maketitle
 \abstract {We formulate and study  the infinite dimensional  linear programming (LP) problem associated with the deterministic discrete time long-run average criterion optimal control problem. Along with its dual, this LP problem allows one to characterize the optimal value of the optimal control problem.  The novelty of our approach is that we focus  on the general  case wherein the optimal value may depend on the initial condition of the system.}

\section{Introduction and Preliminaries}

\ \ \ \ In this paper, we formulate and study  the infinite dimensional (ID) linear programming (LP) problem associated with the deterministic discrete time optimal control problem with  long-run average cost, in which the optimal value may depend on the initial condition of the system. The paper continues the line of research started in \cite{BG-2018}, where similar issues were dealt with in the context of systems evolving in continuous time. Note that, although ideas behind the consideration of continuous and discrete time cases are similar, results in the discrete time case are stronger and are obtained under weaker assumptions comparatively to their continuous time counterparts presented in \cite{BG-2018} (we discuss relationships between the two groups of results in detail in the conclusions section  at the end of the paper).\footnote{An updated and extended version of this paper has been published in SIAM Journal on Control 
and Optimization, Vol. 57, No 3, pp.1783-1817, DOI. 10.1137/18M1229432}

Allowing one  to use the convex duality theory and linear programming  based numerical techniques, LP formulations of various classes of  optimal control problems have been studied extensively  in the literature.
For example,   LP formulations of problems of optimal control of stochastic systems evolving in continuous time   have been considered   in
\cite{Redbook,BhaBor,BGQ,F-V,Kurtz,Stockbridge1}. Various aspects of the LP approach to problems of optimization of discrete time stochastic systems (controlled Markov chains) have been discussed in  \cite{Vivek,Jean,HK-1,HK-2}. In the deterministic  setting, the LP approach has been developed/applied
 in \cite{Goreac-Serea,Her-Her-Lasserre,Lass-Trelat,Rubio,Vinter} for systems evolving in continuous time  considered on a finite time interval. The applicability of the LP approach to deterministic continuous and discrete time systems considered on the infinite time horizon has been explored
 in \cite{GQ,GQ-1,GPS-2017,GPS-2018,QS}.\footnote{Infinite time horizon optimal control problems have been traditionally studied with the help of other (not LP related) techniques; see, e.g.,  \cite{Bardi,Cap-1,Cap-2,CHL,GruneSIAM98,GruneJDE98,Z14,Z06a}  and references therein. Note  that the  list of references mentioned above represents only a sample of the available literature and  is not even  close to being exhaustive. }



 Note that, while the form and the properties of the IDLP problem related to the ergodic case (that is, the case when the  optimal value is independent of the initial conditions) have been well understood, the linear programming formulation of the long-run average optimal control problem in the non-ergodic case has not been discussed much in the literature. In fact,  a justification of counterparts of  LP formulations for reducible finite state Markov chains, as in, e.g., \cite{HK-1} and \cite{HK-2}, presents a significant mathematical challenge. First steps to address this challenge have been made in  \cite{BG-2018}, and (as mentioned above) the present paper is a continuation of this work.

Everywhere in what follows, we will be dealing with the discrete time controlled dynamical system
\begin{equation}\label{A1}
\begin{aligned}
&y(t+1)=f(y(t),u(t)), \; t=0,1,\dots\, \\
&y(0)=y_0,\\
&y(t)\in Y,\\
&u(t)\in U(y(t)).
\end{aligned}
\end{equation}
Here  $Y$ is a given nonempty compact subset of $\reals^m$, $\ U(\cdot):\,Y\leadsto U_0$ is an upper semicontinuous compact-valued mapping to a given compact metric space $U_0$,
$\ f(\cdot, \cdot ):\,\reals^m\times U_0\to \reals^m$ is a continuous function.

It can be observed that the last two constraints of \eqref{A1} can be rewritten as one:
$$
u(t)\in A(y(t)),
$$
where the map $\ A(\cdot):\,Y\leadsto U_0$ is defined by the equation
\begin{equation*}
\begin{aligned}
A(y):=\{u\in U(y)|\, f(y,u)\in Y\} \ \ \ \forall y\in Y.
\end{aligned}
\end{equation*}
The map $A(\cdot)$ is upper semicontinuous and its graph $G$,
$$
G:={\rm graph}\,A=\{(y,u)|\,y\in Y,\,u\in U(y),\,f(y,u)\in Y\},
$$
is a compact subset of $Y\times U_0$.

A control $u(\cdot)$ and the pair $(y(\cdot),u(\cdot))$ will be called an admissible control and an admissible process, respectively, if the relationships \eqref{A1} are satisfied.
The set of admissible controls will be denoted $\U(y_0)$ or $\U_T(y_0)$, depending on whether the problem is considered on the infinite time horizon or on a finite time sequence $t\in \{0,\dots,T-1\}$.

Everywhere in the paper, it is assumed that\\

A1. {\em The set $A(y)$ is not empty for any $y\in Y$.}\\

This assumption implies that the sets  $\U_T(y_0)$ (with $T$ being an arbitrary positive interger) and the set $\U(y_0)$ are not empty for any $y_0\in Y$. That is, there exists at least one admissible control for any initial condition (systems that satisfy such a property are called {\it viable}; see \cite{Aub}).

On  the trajectories of \eqref{A1}, we consider the following optimal control problems:
\begin{equation}\label{A112-1}
{1\o T} \min_{u(\cdot)\in \U_T(y_0)}\sum_{t=0}^{T-1} k(y(t),u(t))=:V_T(y_0),
\end{equation}
\begin{equation}\label{A112-2}
(1-\a) \min_{u(\cdot)\in \U(y_0)}\sum_{t=0}^{\infty}\a^t k(y(t),u(t))=:h_{\a}(y_0),
\end{equation}
where $k:\,\reals^m\times U_0\to \reals^m$ is a continuous function and $\a\in (0,1)$ is a discount factor. Note that, under Assumption A1, the minima in (\ref{A112-1}) and (\ref{A112-2}) are achieved and the optimal value functions $V_T(\cdot) $, $h_{\a}(\cdot)$ are lower  semicontinuous (see, e.g., Propositions 1-3 and Corollary 1 in \cite{GPS-2017}).

An extensive literature is devoted to matters related to the existence  and equality of the limits  $\ \lim_{T\rightarrow\infty}V_T(y_0)$ and $\ \lim_{\a\up 1}h^{\a}(y_0)$. The ergodic case, when these limits are constants (that is, when they do not depend on the initial condition $y_0$),  was studied, for example,  in \cite{Arisawa-3,Redbook,Bardi,GQ} (see also references therein). Results  for the non-ergodic case were obtained in  \cite{BQR-2015,GruneSIAM98,GruneJDE98,Khlopin,Sorin92,OV-2012,QR-2012}. In particular, it was results of  \cite{BQR-2015} that were instrumental for obtaining the IDLP representation for the aforementioned limits for systems evolving in continuous time in \cite{BG-2018}. Some ideas from \cite{BQR-2015}  are used in this paper too.

The paper is organized as follows. In the remainder of this introductory section, we give some definitions and state some earlier results that are used further in the text. In Section \ref{Sec-Bound-from-below}, we introduce an IDLP problem and its dual, the optimal value of the latter giving a lower bound for $\ \liminf_{T\rightarrow\infty}V_T(y_0)$ and $\ \liminf_{\a\up 1}h_{\a}(y_0)$ (see Proposition \ref{Prop-imporant-1}). In Section \ref{Sec-Up-Bound}, we establish (see Theorem \ref{Th-upper-bound}) that  $\ \limsup_{T\rightarrow\infty}V_T(y_0)$ and $\ \limsup_{\a\up 1}h_{\a}(y_0)$ are bounded from above by the optimal value of the IDLP problem introduced in Section \ref{Sec-Bound-from-below} provided that the  value functions $V_T(\cdot) $, $h_{\a}(\cdot)$ are continuous. Note that the proof of Theorem \ref{Th-upper-bound} is based on a lemma that extends some results of \cite{BQR-2015} to the discrete time case (see Lemma \ref{Lem-upper-bound}). A direct corollary from the above mentioned results is Proposition \ref{Prop-lims-exist} of Section \ref{Sec-equality} stating that the limits $\ \lim_{T\rightarrow\infty}V_T(y_0)$ and $\ \lim_{\a\up 1}h^{\a}(y_0)$ exist and are equal to the optimal value of the IDLP problem if there is no duality gap. The main result of
Section \ref{Sec-equality}  is   Theorem \ref{ThN1} establishing  that, if the pointwise limits $\ \lim_{T\rightarrow\infty}V_T(y_0)$ and $\ \lim_{\a\up 1}h^{\a}(y_0)$ exist and are continuous, then they are equal to the optimal value of the dual problem.  Also in this section, we use the optimal solution of the dual IDLP problem to state sufficient and necessary optimality conditions for the long-run average optimal control problem (see Propostions \ref{Prop-optim-cond-suf} and \ref{Prop-optim-cond-nec}), these optimality conditions are illustrated with an elementary \lq\lq toy example". In Section \ref{Sec-Appendix}, we establish some auxiliary results used in the proofs of the previous sections  and in Section \ref{Sec-COnclusions}, we present some  conclusions summarizing results obtained  and comparing them with results of \cite{BG-2018}.

We conclude this section with the introduction of notations and results that are used in the sequel.
Let $(y(\cdot), u(\cdot))$ be an admissible process.
A probability measure
$\g_{(y(\cdot), u(\cdot)),S}$ is called the {\em  occupational measure} generated by the process $(y(\cdot), u(\cdot))$ over the time sequence $\{0,1,...,S-1 \}$ if, for any Borel set $Q\subset G$,
\begin{equation*}
\g_{(y(\cdot), u(\cdot)),S}(Q)={1\o S} \sum_{t=0}^{S-1} 1_Q(y(t),u(t)).
\end{equation*}
A probability measure $\g^{\a}_{(y(\cdot),u(\cdot))}$ is called the {
\em discounted occupational measure} generated by the process $(y(\cdot), u(\cdot))$ if, for any Borel set $Q\subset G$,
\begin{equation}\label{E6}
\g^{\a}_{(y(\cdot),u(\cdot))}(Q)=(1-\a) \sum_{t=0}^{\infty} \a^t 1_Q(y(t),u(t)),
\end{equation}
where $1_Q(\cdot)$ is the indicator function of $Q$.

It can be shown that, if $\g_{(y(\cdot), u(\cdot)),S}$ is the occupational measure generated by the process $(y(\cdot), u(\cdot))$ over the time sequence $\{0,1,...,S-1 \}$, then
\begin{equation}\label{G88}
\int_{G} q(y,u) \g_{(y(\cdot),u(\cdot)),S}(dy,du)={1\o S}\sum_{t=0}^{S-1} q(y(t),u(t))
\end{equation}
for any Borel measurable function $q$ on $G$. Also, it can be shown  that if $\g^{\a}_{(y(\cdot),u(\cdot))}$ is the discounted occupational measure generated by the process $(y(\cdot), u(\cdot))$, then
\begin{equation}\label{G8}
\int_{G} q(y,u) \g^{\a}_{(y(\cdot),u(\cdot))}(dy,du)=(1-\a)\sum_{t=0}^{\infty} \a^t q(y(t),u(t))
\end{equation}
for any Borel measurable function $q$ on $G$.

Let us introduce the following notations for the sets of occupational measures:
\begin{equation}\label{union-1}
\G_T(y_0):=\bigcup_{u(\cdot)\in \U_T(y_0)}\{\g_{(y(\cdot),u(\cdot)),T}\},\quad \ \ \ \ \G_T:=\bigcup_{y_0\in Y}\{\G_T(y_0)\},
\end{equation}
\begin{equation}\label{union-2}
 \T_{\a}(y_0):=\bigcup_{_{u(\cdot)\in \U(y_0)}}\{\g^{\a}_{(y(\cdot),u(\cdot))}\},\quad \ \ \ \ \  \T_{\a}:=\bigcup_{y_0\in Y}\{\T_{\a}(y_0)\}.
\end{equation}
Note that, due to \eqref{G88} and \eqref{G8}, problems \eqref{A112-1} and (\ref{A112-2}) can be rewritten in the form
\begin{equation}\label{A3-T}
\min_{\g\in \G_T(y_0)} \int_{G} k(y,u)\g(dy,du)=V_T(y_0)
\end{equation}
and
\begin{equation}\label{A3}
\min_{\g\in  \T_{\a}(y_0)} \int_{G} k(y,u)\g(dy,du)=(1-\a)h_{\a}(y_0),
\end{equation}
respectively.

To describe convergence properties of occupational measures, we introduce the following metric on $\P(G)$ (the space of probability measures defined on Borel subsets of $G$):
$$
\rho(\g',\g''):=\sum_{j=1}^{\infty} {1\o 2^j}\left|\int_G q_j(y,u)\g'(dy,du)-\int_G q_j(y,u)\g''(dy,du)\right|
$$
for $\g',\g''\in \P(G)$, where $q_j(\cdot),\,j=1,2,\dots,$ is a sequence of Lipschitz continuous functions dense in the unit ball of the space of continuous functions $C(G)$ from $G$ to $\reals$.
This metric is consistent with the weak$^*$ convergence topology on $\P(G)$, that is,
a sequence $\g^k\in \P(G)$ converges to $\g\in \P(G)$ in this metric if and only if
$$
\lim_{k\to \infty}\int_G q(y,u)\g^k(dy,du)=\int_G q(y,u)\g(dy,du)
$$
for any $q\in C(G)$.
Using the metric $\rho$, we can define the ``distance" $\rho(\g,\Gamma)$ between $\g\in \P(G)$ and $\Gamma\subset \P(G)$
and the Hausdorff metric $\rho_H(\Gamma_1,\Gamma_2)$ between $\Gamma_1\subset \P(G)$ and $\Gamma_2\subset \P(G)$ as follows:
$$
\rho(\g,\Gamma):=\inf_{\g'\in \Gamma}\rho(\g,\g'),\quad
\rho_H(\Gamma_1,\Gamma_2):=\max\{\sup_{\g\in \Gamma_1}\rho(\g,\Gamma_2),\sup_{\g\in \Gamma_2}\rho(\g,\Gamma_1)\}.
$$
Note that, although, by some abuse of terminology,  we refer to
$\rho_H(\cdot,\cdot)$ as  a metric on the set of subsets of
${\mathcal P} (G)$, it is, in fact, a semi metric on this set
(since $\rho_H(\Gamma_1, \Gamma_2)=0$ implies  $\Gamma_1
= \Gamma_2$ if  $\Gamma_1$ and $\Gamma_2$ are closed, but the equality may not be true if at least one of these sets is not closed).

Let us define the sets $W$ and $W(\a,y_0)$  by the equations:
\begin{equation*}\label{M17}
\begin{aligned}
W:=\{\g\in \P(G)\ |\, \int_{G}(\ph(f(y,u))-\ph(y))
\gm=0\quad \hbox{for all } \ph\in C(Y)\},
\end{aligned}
\end{equation*}
\begin{equation*}
\begin{aligned}
W(\a,y_0)=&\{\g\in \P(G)|\,\\&\int_G (\alpha\ph(f(y,u))-\ph(y)+(1-\a)(\ph(y_0)-\ph(y)))\gm=0\quad\hbox{for all }\ph\in C(Y)\}.
\end{aligned}
\end{equation*}
Note that the  sets $W$ and $W(\a,y_0)$  are convex and compact in the  topology specified above.
The following equalities establish relationships between these sets and the occupational measures sets introduced earlier (see  Theorem 5.4 in \cite{GPS-2017}):
\begin{equation}\label{convergence-to-W}
\lim_{T\to \infty}\rho_H (\bar{\rm co}\ \G_T,W)=\lim_{\a\up 1}\rho_H (\bar{\rm co}\ \T_{\a},W)=0.
\end{equation}
Also (see Corollary 2 in \cite{GPS-2017}),
\begin{equation}\label{convergence-to-W-dis}
\bar{\rm co}\ \T_{\a}(y_0)=W(\a,y_0) \ \ \ \forall \ \a\in (0,1).
\end{equation}
Here and in what follows, $\bar{\rm co} $ stands for the closed convex hull of the corresponding set.

\section{Estimates of the Limit Optimal Value Functions from Below}\label{Sec-Bound-from-below}

Consider the IDLP problem
\begin{equation}\label{BB1}
\inf_{(\g,\xi)\in \O(y_0)} \int_G k(y,u)\gm=: k^*(y_0),
\end{equation}
where
\begin{equation}\label{eq-Omega}
\begin{aligned}
&\O(y_0):=\{(\g,\xi)\in \P(G)\times \Mp(G)|\, \g\in W,\,\\
&\int_{G}(\ph(y_0)-\ph(y))\gm+\int_G(\ph(f(y,u))-\ph(y))\xm=0
\quad \hbox{for all } \ph\in C(Y)\},
\end{aligned}
\end{equation}
with $\mathcal{M}_+ (G) $ standing for the space of nonnegative measures defined on Borel subsets of $G$.
Also consider the problem
\begin{equation}\label{BB8}
\sup_{(\mu,\psi,\eta)\in \D} \mu=:d^*(y_0),
\end{equation}
where $\D(y_0)$ is the set of triplets $(\mu,\psi(\cdot),\eta(\cdot))\in \reals\times C(Y)\times C(Y)$ that for all $(y,u)\in G$ satisfy the inequalities
\begin{equation}\label{BB7}
\begin{aligned}
&k(y,u)+(\psi(y_0)-\psi(y))+\eta(f(y,u))-\eta(y)-\mu\ge 0,\\
&\psi(f(y,u))-\psi(y)\ge 0.
\end{aligned}
\end{equation}
Note that the optimal value of problem (\ref{BB8})  can be equivalently represented as
\begin{equation}\label{CC9}
d^*(y_0)=\sup_{\psi,\eta}\min_{(y,u)\in G}\{k(y,u)+(\psi(y_0)-\psi(y))+\eta(f(y,u))-\eta(y)\},
\end{equation}
where $\psi$ and $\eta$ are continuous functions, and $\psi$ satisfies the second inequality in \eqref{BB7}.
The optimal values of \eqref{BB8}  and \eqref{BB1} are related by the inequality
\begin{equation}\label{CC9-11}
 d^*(y_0)\leq k^*(y_0)
\end{equation}
(see Lemma \ref{L-equivalince} in Section \ref{Sub-Sec-Auxiliary-Res}). Problem \eqref{BB8} is, in fact, dual with respect to \eqref{BB1}, with
 (\ref{CC9-11}) being a part of the duality relationships (see more details in Section \ref{Sub-Sec-Auxiliary-Res}).


As can be readily seen, problem \eqref{BB1} can be equivalently written as
\begin{equation}\label{CC11}
\inf_{\g\in W_1(y_0)} \int_G k(y,u)\gm=k^*(y_0),
\end{equation}
where
\begin{equation*}
\begin{aligned}
&W_1(y_0)=\{\g\in W|\, \hbox{ there exists } \xi\in\M_{+}(G)\;\hbox{such that }(\g,\xi)\in \O(y_0)\}=\\
&\{\g\in W|\, \hbox{ there exists }\xi\in \Mp(G)\quad \hbox{such that }\\
&\int_G(\ph(y)-\ph(y_0))\gm=\int_G(\ph(f(y,u))-\ph(y))\xm\ \ \ \ \forall \ \ph\in C(Y)\}.
\end{aligned}
\end{equation*}
Along with \eqref{CC11}, consider the problem
\begin{equation}\label{CC14}
\min_{\g\in W_2(y_0)} \int_G k(y,u)\gm,
\end{equation}
where
\begin{equation*}
\begin{aligned}
&W_2(y_0)=\{\g\in W|\, \hbox{ there exists a sequence }\xi_i\in \Mp(G),\,i=1,2,\dots,\quad \hbox{such that }\\
&\int_G(\ph(y)-\ph(y_0))\gm=\lim_{i\to \infty}\int_G(\ph(f(y,u))-\ph(y))\xi_i(dy,du)\ \ \ \forall \ \ph\in C(Y)\}.
\end{aligned}
\end{equation*}
It is easy to see that both sets $W_1(y_0)$ and $W_2(y_0)$ are convex, set $W_2(y_0)$ is closed (and, therefore, compact), and
$$
{\rm cl}\, W_1(y_0)\subset W_2(y_0).
$$

\begin{Lemma}\label{L1}
The following inclusions are true:
\begin{equation}\label{CC3}
\limsup_{T\to \infty} \G_T(y_0)\subset W_2(y_0) \quad\hbox{and }\limsup_{\a\up 1} \T^{\a}(y_0)\subset W_2(y_0).
\end{equation}
This implies, in particular, that the set $W_2(y_0)$ is not empty.
\end{Lemma}

{\bf Proof.} Note first that since the sets $\G_T(y_0)$ and $\T^{\a}(y_0)$ are not empty for all admissible $T$ and $\a$, so are the sets $\disp\limsup_{T\to \infty} \G_T(y_0)$ and $\disp\limsup_{\a\up 1} \T^{\a}(y_0)$.
Note also that from (\ref{convergence-to-W}) it follows that
\begin{equation}\label{CC12}
\limsup_{T\to \infty} \G_T(y_0)\subset W \quad\hbox{and }\limsup_{\a\up 1} \T^{\a}(y_0)\subset W.
\end{equation}
Let $\disp\g\in \limsup_{T\to \infty} \G_T(y_0)$. Then there exist sequences $T_i\to\infty$ and $\g_i\in \G_{T_i}(y_0)$ such that $\g_i\to \g$ as $i\to \infty$. Let $u_i(\cdot)\in \U_{T_i}(y_0)$ be the control generating $\g_i$ and $y_i(\cdot)$ be the corresponding trajectory.  For any $\ph\in C(Y)$ we have
\begin{equation}\label{BB2}
\begin{aligned}
&\int_G (\ph(y)-\ph(y_0))\,\g_i(dy,du)={1\o T_i}\sum_{t=0}^{T_i-1}(\ph(y_i(t))-\ph(y_0))\\
&={1\o T_i}\sum_{t=0}^{T_i-1}\sum_{s=0}^{t-1}(\ph(y_i(s+1))-\ph(y_i(s)))
={1\o T_i}\sum_{t=0}^{T_i-1}\sum_{s=0}^{t-1}(\ph(f(y_i(s),u_i(s)))-\ph(y_i(s))).
\end{aligned}
\end{equation}
Define the functional $\zeta_i\in C^*(G)$ (here and in what follows, $C^*(G)$ stands for the space of continuous linear functionals on $C(G)$) by the equation
$$
\langle \zeta_i,q\rangle ={1\o T_i}\sum_{t=0}^{T_i-1}\sum_{s=0}^{t-1}q(y_i(s),u_i(s)) \quad \hbox{for all }q\in C(G).
$$
Due to Riesz representation theorem (see, e.g., Theorem 4.3.9, p.\ 181 in \cite{Ash}), there exists $\xi_i\in \Mp(G)$ such that
$$
\langle \zeta_i,q\rangle =\int_G q(y,u)\xi_i(dy,du) \quad\hbox{for all }q\in C(G).
$$
Then \eqref{BB2} can be written as
\begin{equation*}
\begin{aligned}
\int_G (\ph(y)-\ph(y_0))\,\g_i(dy,du)=\langle \zeta_i,\ph(f(y,u))-\ph(y)\rangle=\int_G (\ph(f(y,u))-\ph(y))\,\xi_i(dy,du).
\end{aligned}
\end{equation*}
Passing to the limit, we obtain
\begin{equation*}
\begin{aligned}
\int_G (\ph(y)-\ph(y_0))\,\g(dy,du)=\lim_{i\to \infty}\int_G (\ph(f(y,u))-\ph(y))\,\xi_i(dy,du).
\end{aligned}
\end{equation*}
Since $\g\in W$ (due to  \eqref{CC12}), the latter equality implies that $\g\in W_2(y_0).$ Thus, the first inclusion
in \eqref{CC3} is proved.

Let us prove the second inclusion. By (\ref{convergence-to-W-dis}),
 to prove the second inclusion in \eqref{CC3}, it is sufficient to prove that
$$
\limsup_{\a\up 1}W(\a,y_0)\subset W_2(y_0).
$$
Note that from (\ref{convergence-to-W}) and (\ref{convergence-to-W-dis}) it follows that
\begin{equation*}
\limsup_{\a\up 1} W(\a,y_0)\subset W.
\end{equation*}
Take $\g \in \limsup_{\a\up 1}W(\a,y_0)$. There exist sequences $\a_i\up 1$ and $\g_i\in W(\a_i,y_0)$ such that $\g_i\to \g$ as $i\to \infty$. Since $\g_i\in W(\a_i,y_0)$, we have
\begin{eqnarray}\label{CC4}
\int_G (\ph(y)-\ph(y_0))\,\g_i(dy,du) &=& {1\o 1-\a_i}\int_G (\ph(f(y,u))-\ph(y))\,\g_i(dy,du) \nonumber \\
&=& \int_G (\ph(f(y,u))-\ph(y))\,\xi_i(dy,du),
\end{eqnarray}
where $\xi_i=\g_i/(1-\a_i)$.
Passing to the limit as $i\to \infty$ we obtain
\begin{equation*}
\int_G (\ph(y)-\ph(y_0))\,\g(dy,du)=\lim_{i\to \infty}\int_G (\ph(f(y,u))-\ph(y))\,\xi_i(dy,du).
\end{equation*}
Since $\g\in W$, the second inclusion in \eqref{CC3} is proved.\hf

\bigskip

The next lemma establishes a relation between the optimal values in problems \eqref{BB8} and \eqref{CC14}.

\begin{Lemma}\label{L2} The optimal value in problems \eqref{BB8}  and \eqref{CC14} are equal, that is,
$$
d^*(y_0)=\min_{\g\in W_2(y_0)}\int_G k(y,u)\,\gm.
$$
\end{Lemma}

{\bf Proof.} The proof of the lemma is given in Section \ref{Sub-Sec-Auxiliary-Res}. \hf


\begin{Proposition}\label{Prop-imporant-1} The lower limits of the optimal value functions in problems \eqref{A112-1} and \eqref{A112-2} are bounded from below by the optimal value of \eqref{BB8}, that is,
\begin{equation}\label{BB3}
\begin{aligned}
&\liminf_{T\to \infty} V_T(y_0)\ge d^*(y_0),\\
&\liminf_{\a\up 1} h_{\a}(y_0)\ge d^*(y_0).\\
\end{aligned}
\end{equation}
\end{Proposition}

{\bf Proof.} This proposition follows from Lemmas \ref{L1} and \ref{L2}, and from the  fact that the equalities
\begin{equation*}
\begin{aligned}
&\liminf_{T\to \infty} V_T(y_0)=\inf\left\{\int_G k(y,u)\gm,\,\g\in\limsup_{T\to \infty} \G_T(y_0)\right\},\\
&\liminf_{\a\up 1} h_{\a}(y_0)=\inf\left\{\int_G k(y,u)\gm,\,\g\in\limsup_{\a\up 1} \T^{\a}(y_0)\right\}
\end{aligned}
\end{equation*}
are valid. \hf

\medskip

Let $\mathcal{T} $ be a positive integer and let $(y_{\mathcal{T}}(\cdot), u_{\mathcal{T}}(\cdot)) $ be  a $\mathcal{T} $-periodic admissible process. This process will be referred to as {\it finite time (FT) reachable from $y_0$}  if there exist an integer  $\bar t \geq 0 $ and a control $u(\cdot)\in \U_{\bar t}(y_0)$ such that the  solution $y(t)= y(t,y_0,u) $ of (\ref{A1})
obtained with this control satisfies the equality $y(\bar t) = y_{\mathcal{T}}(0) $.

 Consider the optimal control problem
\begin{equation}\label{e-main-8-1}
\inf_{\mathcal{T},\left(y_{\mathcal{T}}(\cdot),u_{\mathcal{T}}(\cdot)\right)  }\left\{\frac{1}{\mathcal{T}}\sum_{t=0}^{\mathcal{T}-1} k(y_{\mathcal{T}}(t),u_{\mathcal{T}}(t)) \right\} :=V_{per}(y_0),
\end{equation}
where ${\rm inf}$ is over all integer $\mathcal{T}>0$ and over all $\mathcal{T}$-periodic pairs $(y_{\mathcal{T}}(\cdot), u_{\mathcal{T}}(\cdot)) $ that are FT reachable from $y_0$. Similarly to (\ref{A3-T}), this problem can be reformulated in terms of occupational measures
\begin{equation}\label{e:occup-meas-def-eq-per}
\inf_{\gamma\in \Gamma_{per}(y_0) }\int_{Y\times U}k(y,u)\gamma(dy,du)=V_{per}(y_0),
\end{equation}
where $\Gamma_{per}(y_0)$ is the set of occupational measures generated by all FT reachable from $y_0 $-admissible  periodic pairs. Note that
\begin{equation}\label{e:occup-meas-def-eq-per-11}
\Gamma_{per}(y_0)\subset \limsup_{T\rightarrow\infty}\Gamma_T(y_0)
\end{equation}
and, therefore,
\begin{equation}\label{e:occup-meas-def-eq-per-1}
V_{per}(y_0)\geq \liminf_{T\rightarrow\infty}V_T(y_0).
\end{equation}

\begin{Proposition}\label{Prop-former-2-4}
The following relationships are valid:
\begin{equation}\label{e-main-8-2}
\Gamma_{per}(y_0)\subset W_1(y_0), \ \ \ \ \ \ V_{per}(y_0)\geq k^*(y_0).
\end{equation}
\end{Proposition}

{\bf Proof.} Due to (\ref{CC11}) and (\ref{e:occup-meas-def-eq-per}), it is sufficient to prove only the first relationship. Note that from (\ref{CC12}) and (\ref{e:occup-meas-def-eq-per-11}) it follows that
\begin{equation}\label{e-main-8-2-5}
\Gamma_{per}(y_0)\subset W.
\end{equation}
Take now an arbitrary $\gamma\in \Gamma_{per}(y_0)$. By definition, it means that $\gamma$ is generated by a $\mathcal{T}$-periodic pair $(y_{\mathcal{T}}(\cdot), u_{\mathcal{T}}(\cdot)) $ that is  FT reachable from $y_0 $. That is, for any continuous function $q(y,u) $,
$$
\int_{ G}q(y,u)\gamma(dy,du)= \frac{1}{\mathcal{T}}\sum_{t=0}^{\mathcal{T}-1}q(y_{\mathcal{T}}(t), u_{\mathcal{T}}(t)).
$$
Consequently, for any $\phi\in C(Y) $,
$$
\int_{G}(\phi(y)-\phi(y_0))\gamma(dy,du)= \frac{1}{\mathcal{T}}\sum_{t=0}^{\mathcal{T}-1}(\phi(y_{\mathcal{T}}(t))-\phi(y_0))
$$
\vspace{-0.4cm}
$$
=\frac{1}{\mathcal{T}}\sum_{t=0}^{\mathcal{T}-1}(\phi(y_{\mathcal{T}}(t))-\phi(y_{\mathcal{T}}(0))) + (\phi(y(\bar t)) -\phi(y_0))
$$
\vspace{-0.2cm}
\begin{equation}\label{e-main-4-per-11}
=\frac{1}{\mathcal{T}}\sum_{t=0}^{\mathcal{T}-1}\left(\sum_{s=0}^{t-1}(\phi(y_{\mathcal{T}}(s+1)) - \phi(y_{\mathcal{T}}(s))) \right) +
\sum_{s=0}^{\bar t-1}\left(\phi(y(s+1)) - \phi(y(s)) \right),
\end{equation}
where $y(t)= y(t,y_0,u) $ is a solution of (\ref{A1}) that satisfies the equality $y(\bar t) = y_{\mathcal{T}}(0) $ (the existence
of $\bar t \geq 0 $ and the existence of a control $u(\cdot)\in \U_{\bar t}(y_0)$ that ensure the validity of this equality follows from the fact that $(y_{\mathcal{T}}(\cdot), u_{\mathcal{T}}(\cdot)) $ is  FT reachable from $y_0 $). Since $y_{\mathcal{T}}(s+1) = f(y_{\mathcal{T}}(s),u_{\mathcal{T}}(s))  $ and $y(s+1) = f(y(s),u(s))  $, from (\ref{e-main-4-per-11}) it follows that
$$
\int_{G}(\phi(y)-\phi(y_0))\gamma(dy,du)
$$
\vspace{-0.2cm}
\begin{equation}\label{e-main-4-per}
= \frac{1}{\mathcal{T}}\sum_{t=0}^{\mathcal{T}-1}\left(\sum_{s=0}^{t-1}(\phi((f(y_{\mathcal{T}}(s),u_{\mathcal{T}}(s))) - \phi(y_{\mathcal{T}}(s))) \right)
+ \sum_{s=0}^{\bar t-1}\left(\phi(f(y(s),u(s))) - \phi(y(s)) \right).
\end{equation}
Define $\zeta\in C^*(G)$ by the equation
$$
\langle \zeta ,q\rangle ={1\o \mathcal{T}}\sum_{t=0}^{\mathcal{T} -1}\sum_{s=0}^{t-1}q(y_\mathcal{T} (s),u_\mathcal{T} (s))
+ \sum_{s=0}^{\bar t-1}q(y (s),u (s)) \ \ \ \ \forall \ q\in C(G).
$$
Due to Riesz representation theorem, there exists $\xi \in \Mp(G)$ such that
$$
\langle \zeta,q\rangle =\int_G q(y,u)\xi (dy,du)  \ \ \ \ \forall \  q\in C(G).
$$
Therefore, (\ref{e-main-4-per}) can be rewritten as
$$
\int_{G}(\phi(y)-\phi(y_0))\gamma(dy,du)= \langle \zeta ,\phi(f(y,u)) - \phi (y)\rangle = \int_{G}(\ph(f(y,u))-\ph(y))\,\xi (dy,du).
$$
Since $\gamma\in W $ (by (\ref{e-main-8-2-5})), the latter implies that $\gamma\in W_1(y_0)$. Thus, the first relationship in
(\ref{e-main-8-2}) is established. \hf

\begin{Corollary}\label{Corollary-former-2-5}
If
\begin{equation}\label{e-main-8-3}
 V_{per}(y_0) = \liminf_{T\rightarrow\infty}V_T(y_0),
\end{equation}
then
$$
\liminf_{T\rightarrow\infty}V_T(y_0)\geq k^*(y_0).
$$
\end{Corollary}

\section{Estimates of the Limit Optimal Value Functions from Above}\label{Sec-Up-Bound}
\begin{Theorem}\label{Th-upper-bound}
{\rm (a)} Let $V_T(\cdot)$ be continuous on $Y$ for all natural $T$. Then
\begin{equation}\label{ub-1}
\limsup_{T\rightarrow\infty} V_{T}(y_0)\leq k^*(y_0) \ \ \forall \ y_0\in Y,
\end{equation}
{\rm (b)} Let $h_{\a}(\cdot)$ be continuous on $Y$ for all $\a\in (0,1)$. Then
\begin{equation}\label{ub-2}
\limsup_{\a\up 1} h_{\a}(y_0)\leq k^*(y_0) \ \ \forall \ y_0\in Y.
\end{equation}
\end{Theorem}
Proof of the theorem is based on the following lemma.

\begin{Lemma}\label{Lem-upper-bound}
For any natural $T$,
\begin{equation}\label{e-before-lim-1}
\int_{G} V_T(y)\,\g(dy,du)\le   \int_{G} k(y,u)\,\g(dy,du)\ \ \ \forall \ \gamma\in W.
\end{equation}
Also, for any $\a\in (0,1) $,
\begin{equation}\label{e-before-lim-2}
\int_{G} h_{\a}(y)\,\g(dy,du)\le   \int_{G} k(y,u)\,\g(dy,du) \ \ \ \forall \ \gamma\in W.
\end{equation}
\end{Lemma}
The proof of the lemma is given at the end of the section.\\

{\bf Proof of Theorem \ref{Th-upper-bound}.}\\

{\it Proof of {\rm (a)}.} Let us fix an arbitrary natural  $T$ and let us consider the following IDLP problem
\begin{equation}\label{BB21-1}
\sup_{(\psi,\eta)\in Q(T)} \psi(y_0)=:  d^*(T,y_0),
\end{equation}
where
$Q(T)$ is the set of pairs $(\psi(\cdot),\eta(\cdot))\in C(Y)\times C(Y)$ that satisfy the inequalities
\begin{equation}\label{BB25-1}
\begin{aligned}
&k(y,u)-\psi(y)+\eta(f(y,u))-\eta(y)\ge 0,\\
&\psi(f(y,u))-\psi(y)\ge -\frac{2M}{T}\ \ \ \ \forall \ (y,u)\in G,
\end{aligned}
\end{equation}
with
\begin{equation}\label{e-M}
M:= \max_{(y,u)\in Y\times U_0}|k(y,u)|.
\end{equation}
Let us show that, for an arbitrary small  $\ve > 0 $, there exists  a function $\eta_{T,\ve}(\cdot)\in C(Y) $ such that
\begin{equation}\label{BB21-1-1}
\left(\psi_{T,\ve}(\cdot), \eta_{T,\ve}(\cdot)\right)\in Q(T), \ \ \ {\rm where} \ \ \psi_{T,\ve}(\cdot):=V_T(\cdot)-\ve.
\end{equation}
Note that, if the inclusion above is established,  it would imply that
\begin{equation}\label{BB24-1}
V_T(y_0)-\ve\le  d^*(T, y_0).
\end{equation}
Let us  first verify that there exists $\eta_{T,\ve}(\cdot)\in C(Y) $ such that the pair $(\psi_{T,\ve}(\cdot), \eta_{T,\ve}(\cdot))$ satisfies the first inequality in (\ref{BB25-1}). To this end, note
that the inequality (\ref{e-before-lim-1}) is equivalent to the inequality
\begin{equation*}
\int_{G}(k(y,u)-V_T(y))\,\gm\ge 0\quad \hbox{for all }\g\in W,
\end{equation*}
which, in turn, is equivalent to
\begin{equation}\label{CC15-1}
\min_{\g\in W}\int_{G}(k(y,u)-V_T(y))\,\gm\ge 0.
\end{equation}
The problem on the left hand side of (\ref{CC15-1}), i.e.,
\begin{equation}\label{CC8-1}
\min_{\g\in W}\int_{G}(k(y,u)-V_T(y))\,\gm,
\end{equation}
is an IDLP problem, its dual being
\begin{equation}\label{BB12-1}
\sup_{\eta\in C(Y)}\inf_{(y,u)\in G}\{k(y,u)-V_T(y)+\eta(f(y,u))-\eta(y)\}.
\end{equation}
The optimal values of \eqref{CC8-1} and \eqref{BB12-1} are equal\ (see Proposition 6 in \cite{GPS-2017}). Therefore, \eqref{CC15-1} is equivalent to
\begin{equation}\label{BB15-1}
\sup_{\eta\in C(Y)}\inf_{(y,u)\in G}\{k(y,u)-V_T(y)+\eta(f(y,u))-\eta(y)\}\ge 0.
\end{equation}
From \eqref{BB15-1} it follows that, for any $\ve>0$, there exists a function $\eta_{T,\ve}(\cdot)\in C(Y)$ such that
\begin{equation}\label{e-eps-feasib-1}
k(y,u)-V_T(y)+\eta_{T,\ve}(f(y,u))-\eta_{T,\ve}(y)\ge -\ve\quad \hbox{for all }(y,u)\in G.
\end{equation}
The latter implies that that the pair $(\psi_{T,\ve}(\cdot), \eta_{T,\ve}(\cdot))$, where $\psi_{T,\ve}(\cdot):= V_T(\cdot)-\ve $, satisfies the first inequality in (\ref{BB25-1}).

Let us now  verify that the function $\psi_{T,\ve}(\cdot)= V_T(\cdot)-\ve$ satisfies the second inequality in (\ref{BB25-1}).
From the dynamic programming principle applied to problem (\ref{A112-1}), it follows that,  for any $T\geq 1$,
\begin{equation}\label{e-before-lim-1-1}
TV_T(y)\le k(y,u)+(T-1)V_{T-1}(f(y,u)) \ \ \ \forall\ (y,u)\in G.
\end{equation}
Also, as can be readily seen,
\begin{equation}\label{e-before-lim-1-1-1}
(T-1)V_{T-1}(y)\le T V_{T }(y) + M \ \ \ \forall\ y\in Y.
\end{equation}
By (\ref{e-before-lim-1-1}) and (\ref{e-before-lim-1-1-1}),
$$
TV_T(y)\le k(y,u) +  TV_{T }(f(y,u)) + M \leq  TV_{T }(f(y,u)) +  2M.
$$
Consequently,
$$
V_T(y)\leq V_{T }(f(y,u)) +  \frac{2M}{T}\ \ \ \Rightarrow \ \ \ \psi_{T,\ve}(y)\leq \psi_{T,\ve}(f(y,u)) +  \frac{2M}{T}
$$
Thus, $\psi_{T,\ve}(\cdot)= V_T(\cdot)-\ve$ satisfies the second inequality in (\ref{BB25-1}). Hence,
(\ref{BB21-1-1}) is valid and, consequently, (\ref{BB24-1}) is valid too.

By Lemma \ref{L-equivalince} of Section \ref{LP-Duality-results},
\begin{equation}\label{e-inequality-lemma-1}
 d^*(T, y_0)\leq k^*(T, y_0),
\end{equation}
where
\begin{equation}\label{BB1-T}
 k^*(T, y_0)=\inf_{(\g,\xi)\in \O(y_0)}\left\{ \int_G k(y,u)\gm  + \frac{2M}{T}\int_G\xi(dy,du)\right\}.
\end{equation}
(Note that, to adjust the notations used above and those used in Lemma \ref{L-equivalince}, one should write $ d^*(T, y_0)$ and $ k^*(T, y_0)$ as $d^*(\theta_T, y_0) $ and
$ k^*(\theta_T, y_0)$, where $\theta_T= \frac{2M}{T} $.)

From (\ref{BB24-1}) and (\ref{e-inequality-lemma-1}) it follows that
$\
V_T(y_0)-\ve \le k^*(T, y_0),
$
which implies that
\begin{equation}\label{BB1-T-1}
V_T(y_0)\le k^*(T, y_0)
\end{equation}
since $\ve > 0 $ is arbitrary small. Due to (\ref{BB1-T-1}), to prove (\ref{ub-1}), it is sufficient to establish that
\begin{equation}\label{BB1-T-2}
\lim_{T\rightarrow\infty}k^*(T, y_0)= k^*(y_0).
\end{equation}
One can readily see that $k^*(T, y_0) $ is a decreasing function of $T$ and that $\ k^*(T, y_0) \geq k^*( y_0)   $ for any $T\geq 1 $.
Hence,
$$
\lim_{T\rightarrow\infty}k^*(T, y_0)\geq k^*(y_0).
$$
Let us now show that the opposite inequality is also valid. Let $\delta >0 $ be arbitrary small
 and let $(\gamma', \xi')\in \Omega(y_0) $
be $\delta$-optimal for (\ref{BB1}).
That is,
$$
 \int_{Y\times U}k(y,u)\gamma'(dy,du)\leq k^*(y_0) + \delta.
$$
Then
$$
k^*(y_0, T)\leq \int_{Y\times U} k(y,u)\gamma'(dy,du)+\frac{2M}{T}\int_{Y\times U} \xi'(dy,du)
$$
$$
\leq  k^*(y_0) + \delta  + \frac{2M}{T}\int_{Y\times U}\xi'(dy,du),
$$
$$
\Rightarrow\ \ \ \ \lim_{T\rightarrow \infty}k^*(y_0, T)\leq k^*(y_0) + \delta \ \ \ \Rightarrow \ \ \ \lim_{T\rightarrow \infty}k^*(y_0, T)\leq k^*(y_0)
$$
($\delta > 0 $ can be arbitrary small). Thus (\ref{BB1-T-2}) is established and statement (a) is proved.\\

{\it Proof of {\rm (b)}}
The proof of {\rm (b)} is very similar to that of {\rm (a)}. We fix an arbitrary   $\a\in (0,1)$ and consider the IDLP problem
\begin{equation}\label{BB21-1-2}
\sup_{(\psi,\eta)\in Q(\a)} \psi(y_0)=:  d^*(\a,y_0),
\end{equation}
where
$Q(\a)$ is the set of pairs $(\psi(\cdot),\eta(\cdot))\in C(Y)\times C(Y)$ that satisfy the inequalities
\begin{equation}\label{BB25-1-2}
\begin{aligned}
&k(y,u)-\psi(y)+\eta(f(y,u))-\eta(y)\ge 0,\\
&\psi(f(y,u))-\psi(y)\ge -2M(1-\a)\ \ \ \ \forall \ (y,u)\in G.
\end{aligned}
\end{equation}
We then show that, for an arbitrary small  $\ve > 0 $, there exists  a function $\eta_{\a,\ve}(\cdot)\in C(Y) $ such that
\begin{equation}\label{BB21-1-3}
\left(\psi_{\a,\ve}(\cdot), \eta_{\a,\ve}(\cdot)\right)\in Q(\a), \ \ \ {\rm where} \ \ \psi_{\a,\ve}(\cdot):=h_{\a}(\cdot)-\ve,
\end{equation}
with the inclusion above  implying that
\begin{equation}\label{BB24-3}
h_{\a}(y_0)-\ve\le  d^*(\a, y_0).
\end{equation}
To verify (\ref{BB21-1-3}), we first show
that there exists $\eta_{\a,\ve}(\cdot)\in C(Y) $ such that the pair $(\psi_{\a,\ve}(\cdot), \eta_{\a,\ve}(\cdot))$ satisfies the first inequality in (\ref{BB25-1-2}). As in the proof of {\rm (a)}, we rewrite
 the inequality (\ref{e-before-lim-2}) in the form
\begin{equation*}
\int_{G}(k(y,u)-h_{\a}(y))\,\gm\ge 0 \ \ \ \ \forall \ \g\in W,
\end{equation*}
which is equivalent to\begin{equation}\label{CC15-3}
\min_{\g\in W}\int_{G}(k(y,u)-h_{\a}(y))\,\gm\ge 0.
\end{equation}
The problem on the left hand side of (\ref{CC15-3}), i.e.,
\begin{equation}\label{CC8-3}
\min_{\g\in W}\int_{G}(k(y,u)-h_{\a}(y))\,\gm,
\end{equation}
is an IDLP problem, the  dual of which is
\begin{equation}\label{BB12-3}
\sup_{\eta\in C(Y)}\inf_{(y,u)\in G}\{k(y,u)-h_{\a}(y)+\eta(f(y,u))-\eta(y)\}.
\end{equation}
The optimal values of \eqref{CC8-3} and \eqref{BB12-3} are equal (Proposition 6 in \cite{GPS-2017}). Therefore, \eqref{CC15-3} is equivalent to
\begin{equation}\label{BB15-1-3}
\sup_{\eta\in C(Y)}\inf_{(y,u)\in G}\{k(y,u)-h_{\a}(y)+\eta(f(y,u))-\eta(y)\}\ge 0.
\end{equation}
From \eqref{BB15-1-3} it follows that, for any $\ve>0$, there exists a function $\eta_{\a,\ve}(\cdot)\in C(Y)$ such that
\begin{equation}\label{e-eps-feasib-1-3}
k(y,u)-h_{\a}(y) +\eta_{\a,\ve}(f(y,u))-\eta_{\a,\ve}(y)\ge -\ve\ \ \ \ \forall \ (y,u)\in G.
\end{equation}
The latter implies that  the pair $(\psi_{\a,\ve}(\cdot), \eta_{\a,\ve}(\cdot))$, where $\psi_{\a,\ve}(\cdot):= h_{\a}(\cdot)-\ve $, satisfies the first inequality in (\ref{BB25-1-2}).

To  verify that the function $\psi_{\a,\ve}(\cdot)= h_{\a}(\cdot)-\ve$ satisfies the second inequality in (\ref{BB25-1-2}), note that
from the dynamic programming principle applied to problem (\ref{A112-2}), it follows that
\begin{equation}\label{e-before-lim-2-1}
h_{\a}(y)\le (1-\a)k(y,u)+\a h_{\a}(f(y,u))\ \ \ \forall \ (y,u)\in G
\end{equation}
(see, e.g., Proposition 4 in \cite{GPS-2017}). The latter implies that
$$
h_{\a}(y)\le h_{\a}(f(y,u)) + (1-\a)(k(y,u) - h_{\a}(f(y,u))),
$$
which, in turn, implies that
\begin{equation}\label{e-before-lim-2-2}
h_{\a}(y)\le  h_{\a}(f(y,u)))+ 2M(1-\a) \ \ \ \forall \ (y,u)\in G
\end{equation}
(since, as can be readily seen, $\ \max_{y\in Y}|h_{\a}(y)|\leq M $).
Thus, $\psi_{\a,\ve}(\cdot)= h_{\a}(\cdot)-\ve$ satisfies the second inequality in (\ref{BB25-1-2}), and, therefore,
(\ref{BB24-3}) is valid too.
Starting from this point, the proof of {\rm (b)} follows exactly the same steps as that of {\rm (a)}.  \hf

\bigskip

{\bf Proof of Lemma \ref{Lem-upper-bound}.}
Let us prove  (\ref{e-before-lim-1}). To this end, let us show first that, for any natural $T$ and $T'$,
\begin{equation}\label{CC20}
\int_{G} V_T(y)\,\g'(dy,du)\le   \int_{G} k(y,u)\,\g'(dy,du)+{M(T-1)\o T'}\quad \ \forall \ \g'\in  \G_{T'}(y_0), \ \forall \ y_0\in Y,
\end{equation}
where $M$ is as in (\ref{e-M}).
Take $y_0\in Y$, $\g'\in \G_{T'}(y_0)$, and let  $u(\cdot)\in \U_{T'}(y_0)$ be a control that generates $\g'$ on $\{0,\dots,T'-1\}$. Extend $u$ from the interval $\{0,\dots,T'-1\}$ to the interval $\{0,\dots,T'+T-1\}$ so that $u\in \U_{T'+T}(y_0)$. Such extension is possible due to viability of $Y$. Let $y(\cdot)$ be the corresponding trajectory. Taking into account that
$\disp V_T(y(s))\le  {1\o T}\sum_{r=0}^{T-1} k(y(r+s),u(r+s)))$ for all $s\in \{0,\dots,T'-1\}$, we obtain
\begin{equation*}
\begin{aligned}
&\int_{G} V_T(y)\,\g'(dy,du)={1\o T'} \sum_{s=0}^{T'-1} V_T(y(s))
\le {1\o T'} \sum_{s=0}^{T'-1}{1\o T} \sum_{r=0}^{T-1} k(y(r+s),u(r+s)))\\
&={1\o T} \sum_{r=0}^{T-1}{1\o T'} \sum_{s=0}^{T'-1} k(y(r+s),u(r+s))\\
&={1\o T} \sum_{r=0}^{T-1}{1\o T'} \sum_{\s=r}^{T'+r-1} k(y(\s),u(\s))
\le {1\o T} \sum_{r=0}^{T-1}{1\o T'} \left(\sum_{\s=0}^{T'-1} k(y(\s),u(\s))+2Mr\right)\\
&= {1\o T} \sum_{r=0}^{T-1}{1\o T'} \sum_{\s=0}^{T'-1} k(y(\s),u(\s))+ {1\o TT'} \sum_{r=0}^{T-1}2Mr\\
&={1\o T} \sum_{r=0}^{T-1}\int_{G} k(y,u)\,\g'(dy,du)+{M(T-1)\o T'}=\int_{G} k(y,u)\,\g'(dy,du)+{M(T-1)\o T'}\,.
\end{aligned}
\end{equation*}
Thus the inequality (\ref{CC20}) is established. From this inequality
it follows that
\begin{equation}\label{CC20-2}
\int_{G} V_T(y)\,\g'(dy,du)\le   \int_{G} k(y,u)\,\g'(dy,du)+{M(T-1)\o T'}\quad \ \ \ \forall\g'\in {\rm co}\ \G_{T'}\ ,
\end{equation}
where $\G_{T'} $ is the union of $\G_{T'}(y_0) $ over $y_0\in Y$ (see (\ref{union-1})).
 Take an arbitrary $\gamma\in W $. From (\ref{convergence-to-W}) it follows that there exist sequences $T'_l > 0, \ \g'_l\in \G_{T'_l} $, $\ l=1,2,...,$ such that $T'_l\rightarrow\infty $ and $\g'_l\rightarrow \g $. Passing to the limit along these sequences in (\ref{CC20-2})
and having in mind that
$$
\int_{G} V_T(y)\,\g(dy,du)\leq \liminf_{\g'_l\rightarrow \g }\int_{G} V_T(y)\,\g'_l(dy,du)
$$
(since $V_T (\cdot)  $ is lower  semicontinuous for any $T>0$; see, e.g., Theorem 3.1.5 in \cite{Strook}), one arrives at  inequality (\ref{e-before-lim-1}).

Let us now prove (\ref{e-before-lim-2}).  To this end, let us show first that, for any $\a\in (0,1)$ and any $\a'\in (\alpha,1)$,
\begin{equation}\label{e-corollary-2-1}
\int_{G} h_{\a}(y)\,\g'(dy,du)\le   \frac{1-\alpha}{1-\frac{\alpha}{\alpha'}}\int_{G} k(y,u)\,\g'(dy,du)
\end{equation}
$$
+\left(\frac{1-\alpha}{1-\frac{\alpha}{\alpha'}}-1\right)M\
\ \ \forall \ \g'\in \Theta_{\a'}(y_0), \ \ \forall y_0\in Y.
$$
Take $y_0\in Y$, $\g'\in \T_{\a'}(y_0)$, and let  $u(\cdot)\in \U(y_0)$ be a control that generates $\g'$. Let also $y(\cdot)$ be the trajectory corresponding to $u(\cdot)$. We have
\begin{equation*}
\begin{aligned}
&\int_{G} h_{\a}(y)\,\g'(dy,du)=(1-\a') \sum_{s=0}^{\infty}(\a')^s h_{\a}(y(s))\\
&\le (1-\a') \sum_{s=0}^{\infty}(\a')^s(1-\a) \sum_{r=0}^{\infty} \a^r k(y(r+s),u(r+s))\\
&=(1-\a) \sum_{r=0}^{\infty} \a^r(1-\a')\sum_{s=0}^{\infty} (\a')^s k(y(r+s),u(r+s))\\
&=(1-\a) \sum_{r=0}^{\infty} \a^r(1-\a')(\a')^{-r}\sum_{\s=r}^{\infty} (\a')^{\sigma} k(y(\s),u(\s))\\
&\le (1-\a) \sum_{r=0}^{\infty} \a^r(1-\a')(\a')^{-r}\left(\sum_{\s=0}^{\infty} (\a')^{\s} k(y(\s),u(\s))+ \sum_{\s=0}^{r-1} (\a')^{\s}M\right)\\
&= (1-\a) \sum_{r=0}^{\infty} \a^r (\a')^{-r} \int_{G} k(y,u)\,\g'(dy,du)+ (1-\a) \sum_{r=0}^{\infty} \a^r (\a')^{-r}(1-(\alpha')^{r})M\\
&= \frac{1-\alpha}{1-\frac{\alpha}{\alpha'}}\int_{G} k(y,u)\,\g'(dy,du)
+\left(\frac{1-\alpha}{1-\frac{\alpha}{\alpha'}}-1\right)M.
\end{aligned}
\end{equation*}
From (\ref{e-corollary-2-1}) it follows that
\begin{equation}\label{e-corollary-2-2}
\int_{G} h_{\a}(y)\,\g'(dy,du)\le   \frac{1-\alpha}{1-\frac{\alpha}{\alpha'}}\int_{G} k(y,u)\,\g'(dy,du)
+\left(\frac{1-\alpha}{1-\frac{\alpha}{\alpha'}}-1\right)M\ \ \ \forall \ \g'\in {\rm co}\ \T_{\a'},
\end{equation}
where $\T_{\a'} $ is the union of $\T_{\a'}(y_0) $ over $y_0\in Y$ (see (\ref{union-2})). Take an arbitrary $\gamma\in W $. From (\ref{convergence-to-W}) it follows that there exist sequences $\a'_l \in (0,1), \ \g'_l\in \G_{\a'_l} $, $\ l=1,2,...,$ such that $\a'_l\up 1 $ and $\g'_l\rightarrow \g $. Passing to the limit along these sequences in (\ref{e-corollary-2-2})
and keeping in mind that
$$
\int_{G} h_{\a}(y)\,\g(dy,du)  \leq \liminf_{\g'_l\rightarrow \g}\int_{G} h_{\a}(y)\,\g'_l(dy,du)
$$
(since
$h_{\a} (\cdot)  $ is lower  semicontinuous for any $\a\in(0,1) $; see also Theorem 3.1.5 in \cite{Strook}), one arrives at  inequality (\ref{e-before-lim-2}). \hf

\section{LP Representation for the  Optimal Value and Related Sufficient/Necessary Optimality Conditions}\label{Sec-equality}
The following statement is a direct corollary of Theorem \ref{Th-upper-bound} and Proposition \ref{Prop-imporant-1}.
\begin{Proposition}\label{Prop-lims-exist}
If
\begin{equation}\label{e-strong-duality}
 d^*(y_0)= k^*(y_0),
\end{equation}
then, provided that $V_{T}(\cdot) $ is continuous for any $T>1$, there exists the pointwise limit
\begin{equation}\label{lim-exists-1}
\lim_{T\rightarrow\infty} V_{T}(y_0)= d^*(y_0) \ \ \ \forall \ y_0\in Y.
\end{equation}
Also, provided that $h_{\a}(\cdot) $ is continuous for any $\a\in (0,1) $, there exists the pointwise limit
\begin{equation}\label{lim-exists-2}
\lim_{\a\rightarrow 0} h_{\a}(y_0)= d^*(y_0) \ \ \ \forall \ y_0\in Y.
\end{equation}
\end{Proposition}

Note that a statement about the LP representation of the pointwise limits (\ref{lim-exists-1}) and (\ref{lim-exists-2}) can be established without the  strong duality assumption (\ref{e-strong-duality}) . Namely, the following result is valid.
\begin{Theorem}\label{ThN1}
{\rm (a)} Let the pointwise limit
\begin{equation}\label{chesaro-lim-1-1}
\lim_{T\rightarrow\infty} V_{T}(y_0):= V(y_0) \ \ \ \forall \ y_0\in Y.
\end{equation}
 exist and let the function
 $V (\cdot)  $ be continuous.  Then
\begin{equation}\label{CC10}
 V(y_0)= d^*(y_0) \ \ \forall \ y_0\in Y.
\end{equation}
{\rm (b)} Let the pointwise limit
\begin{equation}\label{abel-lim-3-1}
\lim_{\a\rightarrow 1} h_{\a}(y_0):= h(y_0) \ \ \forall \ y_0\in Y,
\end{equation}
 exist and
the function $h (\cdot)  $ be continuous. Then
\begin{equation}\label{abel-lim-4-1}
 h(y_0)= d^*(y_0) \ \ \ \forall \ y_0\in Y.
\end{equation}
\end{Theorem}
{\bf Proof.} The proof of the theorem is given at the end of this section. \hf

\begin{Remark}\label{Rem-extra}
{\rm If (\ref{chesaro-lim-1-1}) and (\ref{CC10}) are valid, then  the strong duality equality (\ref{e-strong-duality}) is true provided that  condition (\ref{e-main-8-3}) of Corollary \ref{Corollary-former-2-5} is satisfied.  }
\end{Remark}

In the rest of this section, we assume
 that the pointwise limit $\lim_{T\rightarrow\infty}V_T(\cdot)=V(y) $ exists and is continuous, and, therefore, it
 is equal  to the optimal value $d^*(y_0)$ of the dual problem (\ref{BB8}) (by Theorem \ref{ThN1}). That is, (\ref{chesaro-lim-1-1}) and (\ref{CC10}) are valid.

 Consider the optimal control problem
\begin{equation}\label{CC21}
\inf_{u(\cdot)\in \U(y_0)} \liminf_{T\to \infty} {1\o T}\sum_{t=0}^{T-1} k(y(t),u(t))=V(y_0).
\end{equation}
Note that, due to (\ref{chesaro-lim-1-1}),  the optimal value of (\ref{CC21}) is equal to $V(y_0)$ (see Proposition \ref{PN3} in Section \ref{LP-Duality-results}).
Below, we discuss   sufficient and necessary optimality conditions for   problem (\ref{CC21})
stated in terms of an optimal solution of problem \eqref{BB8}.\\

DEFINITION. A pair $(\bar \psi(\cdot), \bar \eta(\cdot))\in C(G)\times C(G) $ will be called {\it an optimal solution of
(\ref{BB8})} if it satisfies the inequalities (compare with (\ref{BB7}))
\begin{equation}\label{BB7-3}
\begin{aligned}
&k(y,u)+(\bar\psi(y_0)-\bar\psi(y))+\bar\eta(f(y,u))-\bar\eta(y)\ge d^*(y_0),\\
&\bar\psi(f(y,u))-\bar\psi(y)\ge 0.
\end{aligned}
\end{equation}

\bigskip


\begin{Proposition}\label{PN2}
{\rm (a)} A  pair  $ (\bar \psi(\cdot), \bar \eta(\cdot)) $ is an optimal solution of (\ref{BB8}) if and only if $\ \bar \psi(\cdot)$ satisfies the second inequality in (\ref{BB7-3}) and
\begin{equation}\label{min}
\begin{aligned}
\min_{(y,u)\in G}\{k(y,u)-\bp(y)+\bar\eta(f(y,u))-\bar\eta(y)\}=V(y_0)-\bp(y_0).
\end{aligned}
\end{equation}
{\rm (b)} If $\bar \eta(\cdot) $ is such that
\begin{equation}\label{min-V}
\begin{aligned}
\min_{(y,u)\in G}\{k(y,u)- V(y)+\bar\eta(f(y,u))-\bar\eta(y)\}=0,
\end{aligned}
\end{equation}
 then the pair
$ ( \bar \psi(\cdot), \bar \eta(\cdot)) $, where $\bar \psi(\cdot) = V(\cdot) $,
is an optimal solution of problem  (\ref{BB8}).

\end{Proposition}
{\bf Proof.} By (\ref{CC9}), the first inequality in (\ref{BB7-3}) is equivalent to the equality
\begin{equation}\label{min-0}
\min_{(y,u)\in G}\{k(y,u)+\bp(y_0)-\bp(y)+\bar\eta(f(y,u))-\bar\eta(y)\}=d^*(y_0).
\end{equation}
Also,
(\ref{min-0}) is equivalent to  (\ref{min}) (due to  (\ref{CC10})). Thus {\rm (a)} is proved.

If $\bar \eta(\cdot) $ is such that (\ref{min-V}) is satisfied, then the pair $ ( \bar \psi(\cdot), \bar \eta(\cdot)) $, where $\bar \psi(\cdot) = V(\cdot) $, satisfies (\ref{min}). Therefore, by {\rm (a)}, this pair is an optimal solution of  (\ref{BB8}). This proves
{\rm (b)}. \hf

\begin{Proposition}\label{Prop-optim-cond-suf}
 Let an optimal solution $(\bar\psi(\cdot),\bar\eta(\cdot))$  of  \eqref{BB8} exist.
Then, for an admissible process $(y(\cdot),u(\cdot))$ to be optimal in \eqref{CC21} it is sufficient that the  equalities
\begin{equation}\label{e-opt-1}
k(y(t),u(t))-\bp(y(t))+\bar\eta(f(y(t),u(t)))-\bar\eta(y(t))=V(y_0)-\bp(y_0),
\end{equation}
\begin{equation}\label{e-opt-2}
\bar\psi(y(t))= \bar\psi(y_0)
\end{equation}
are satisfied for all $t=0,1,...\ $.
\end{Proposition}
{\bf Proof.} From (\ref{e-opt-1}) and (\ref{e-opt-2}) it follows that
$$
k(y(t),u(t))+\bar\eta(f(y(t),u(t)))-\bar\eta(y(t))=V(y_0)
$$
for  all $t=0,1,...\ $.
Therefore, for any $T\geq 1$,
\begin{equation}\label{For-Remark}
\begin{aligned}
&{1\o T}\sum_{t=0}^{T-1}\Big(k(y(t),u(t))+\bar\eta(f(y(t),u(t)))-\bar\eta(y(t))\Big)\\
&={1\o T}\sum_{t=0}^{T-1}\Big(k(y(t),u(t))+\bar\eta(y(t+1))-\bar\eta(y(t))\Big)\\
&={1\o T}\sum_{t=0}^{T-1}k(y(t),u(t))+{1\o T}(\bar\eta(y(T)))-\bar\eta(y(0)))=V(y_0).
\end{aligned}
\end{equation}
Taking into account that
$$
\lim_{T\to \infty}{1\o T}(\bar\eta(y(T))-\bar\eta(y(0)))=0,
$$
 we obtain
$$
\lim_{T\to \infty}{1\o T}\sum_{t=0}^{T-1}k(y(t),u(t))=V(y_0).
$$
That is,
the process $(y(\cdot),u(\cdot))$ is optimal in \eqref{CC21}. \hf

We will now establish that the fulfillment of (\ref{e-opt-1})-(\ref{e-opt-2}) is also a necessary condition of optimality of an admissible process $(y(\cdot),u(\cdot))$ provided that the latter is periodic, that is, there exists a positive integer $T_0$ such that, for any $t=0,1,... $,
\begin{equation}\label{e-period}
(y(t),u(t)) = (y(t+T_0),u(t+T_0)) \ \ \ \ \forall\ t=0,1,... \ .
\end{equation}

\begin{Proposition}\label{Prop-optim-cond-nec}
 Let an optimal solution $(\bar\psi(\cdot),\bar\eta(\cdot))$  of  \eqref{BB8} exist. Then, for an  admissible process
 $(y(\cdot),u(\cdot))$ satisfying the periodicity conditions {\rm (\ref{e-period})} to be optimal in {\rm (\ref{CC21})}, it is necessary that the equalities {\rm (\ref{e-opt-1})-(\ref{e-opt-2})}
 are satisfied for all $t=0,1,... $.
 \end{Proposition}
{\bf Proof.}
Note that the fact that the periodic admissible process is optimal in (\ref{CC21}) means that
\begin{equation}\label{e-period-opt-1}
{1\o T_0}\sum_{t=0}^{T_0-1} k(y(t),u(t)) = V(y_0).
\end{equation}
Note also that from Proposition \ref{PN2} it follows that, for any $t=0,1,...,T_0-1$,
\begin{equation}\label{e-period-opt-2}
k(y(t),u(t))-\bp(y(t))+\bar\eta(f(y(t),u(t)))-\bar\eta(y(t))\geq V(y_0)-\bp(y_0),
\end{equation}
\begin{equation}\label{e-period-opt-3}
\bp(y(t))\geq \bp(y_0)
\end{equation}
From (\ref{e-period-opt-1}) and (\ref{e-period-opt-2}) it follows that
$$
\sum_{t=0}^{T_0-1} (\bp(y_0)-\bp(y(t))) + \sum_{t=0}^{T_0-1}(\bar\eta(f(y(t),u(t)))-\bar\eta(y(t)))\geq 0,
$$
which implies that
\begin{equation}\label{e-period-opt-4}
\sum_{t=0}^{T_0-1} (\bp(y_0)-\bp(y(t)))\geq 0
\end{equation}
due to the fact that
\begin{equation}\label{e-period-opt-5}
\sum_{t=0}^{T_0-1}(\bar\eta(f(y(t),u(t)))-\bar\eta(y(t)))= \sum_{t=0}^{T_0-1}(\bar\eta(y(t+1))-\bar\eta(y(t))) = \bar\eta(y(T_0))- \bar\eta(y_0)=0
\end{equation}
(by (\ref{e-period})). The inequalities (\ref{e-period-opt-3}) and (\ref{e-period-opt-4}) establish the validity of (\ref{e-opt-2}).
In view of (\ref{e-opt-2}), the inequality (\ref{e-period-opt-2}) is equivalent to that
\begin{equation}\label{e-period-opt-6}
k(y(t),u(t))+\bar\eta(f(y(t),u(t)))-\bar\eta(y(t))\geq V(y_0)
\end{equation}
for all $t=0,1,...,T_0-1$.
If the above inequality was strict for at least one $t$, then one would obtain
$$
{1\o T_0} \sum_{t=0}^{T_0-1}( k(y(t),u(t))+\bar\eta(f(y(t),u(t)))-\bar\eta(y(t)))>  V(y_0),
$$
which, by (\ref{e-period-opt-5}), would lead to
$$
{1\o T_0}\sum_{t=0}^{T_0-1} k(y(t),u(t)) >  V(y_0).
$$
The latter contradicts (\ref{e-period-opt-1}).
Hence, (\ref{e-period-opt-6}) is satisfied as equality for all $t=0,1,...T_0-1$. This proves (\ref{e-opt-1}).
\hf

\begin{Remark}\label{Rem-feedback}
{\rm As established by Proposition \ref{Prop-optim-cond-suf}, an admissible  process $(y(\cdot),u(\cdot)) $ is optimal if it satisfies the equalities (\ref{e-opt-1}), (\ref{e-opt-2}). Assuming that these are valid, one may conclude  (due to (\ref{min})) that the equality (\ref{e-opt-1}) is equivalent to
$$
(y(t),u(t))={\rm argmin}_{(y,u)\in G}\{k(y,u)-\bar\psi(y)+\bar\eta(f(y,u))-\bar\eta(y)\}
$$
which leads to
$$
u(t)={\rm argmin}_{u\in A(y)}\{k(y(t),u)+\bar\eta(f(y(t),u))\} \ \ \ \ \forall t=0,1,...\ .
$$
The latter implies that the  feedback control
\begin{equation}\label{e-feedback}
\begin{aligned}
u(y)={\rm argmin}_{u\in A(y)}\{k(y ,u)+\bar\eta(f(y ,u))\}
\end{aligned}
\end{equation}
is optimal in the sense that, being used in  (\ref{A1}), it allows one to obtain the optimal \lq\lq open loop" admissible process   $(y(\cdot),u(\cdot)) $.}
\end{Remark}

Let us illustrate the optimality conditions discussed above with the following \lq\lq toy example".\\

{\bf Example.} Let the dynamics be one-dimensional and be described by the  equation (compare with (\ref{A1}))
$$
y(t+1)= u(t)y(t) \ \ \ \forall \ t=0,1,...\ ,
$$
with
$
\ Y=[-1,1]
$ and with $U(y)=\{-1, 1\}$ (that is, the control can be either equal to $1$ or to $-1$). Consider problem (\ref{A112-1}) with
$\
k(y,u)=y
$. As can be readily understood,  the optimal admissible processes  in this example are as follows. If $y_0\in (0,1]$, then
$$
 \ u(0)=-1, \ \  y(0)=y_0 \ \ \ \ {\rm and} \ \ \ \ u(t)=1, \ \ y(t)=- y_0 \ \ \ \ \forall\ \ \ t\geq 1.
 $$
 If $y_0\in [-1,0)$, then
 $$
  u(t)=1, \ \ y(t)= y_0 \ \ \ \ \forall\ \ \ t\geq 0.
 $$
Also, if $y_0= 0 $, then the system is uncontrollable, and the only admissible trajectory is $\ y(t)=0 \ \forall \ t\geq 0 $. The admissible processes described above are optimal on any time horizon (both finite and infinite), with the optimal value function being defined by the equation
\begin{equation}\label{e-example-11-0}
V_T(y_0) = \frac{1}{T}y_0 - \frac{T-1}{T}y_0 = -y_0 + \frac{2}{T}y_0 \ \  {\rm if} \ \  y_0\in (0,1]\ \ \ \ \ {\rm and} \ \ \ \ V_T(y_0) = y_0 \ \ \ {\rm if} \ \ \ y_0\in [-1,0].
\end{equation}
Thus,
 $V(y )= - |y | $. Note that condition (\ref{e-main-8-3}) of Corollary \ref{Corollary-former-2-5} is satisfied and, therefore, the strong duality equality (\ref{e-strong-duality}) is valid in the given example (see Remark \ref{Rem-extra}).

 Define  the function $\bar{\eta}(\cdot)$ by the equation
\begin{equation}\label{e-example-11}
 \bar{\eta}(y):=\max\{2y , 0\} \ \ \ \forall \ y \in [-1,1].
\end{equation}
One can readily verify that
$$
\min_{u\in \{-1,1 \}}\bar{\eta}(uy)= 0 \ \ \ \forall\ y\in [-1,1],
$$
the latter implying that
$$
y+|y| + \min_{u\in \{-1,1 \}}\bar{\eta}(uy) - \bar\eta(y) = 0 \ \ \ \forall y\in [-1,1].
$$
That is,  $\bar{\eta}(\cdot)$  satisfies  (\ref{min-V}).
Therefore, the pair $(\bar\psi(\cdot),  \bar{\eta}(\cdot) ) $, where $\bar\psi(y)=- |y| $,  is an optimal solution of (\ref{BB8}). The $argmin$ feedback control defined in (\ref{e-feedback})  takes in this case the form
$$
{\rm argmin}_{u\in \{-1,1 \}}\bar{\eta}(uy) = -1  \ \ {\rm if } \ \ y\in (0,1), \ \ \ \ \ {\rm argmin}_{u\in \{-1,1 \}}\bar{\eta}(uy) = 1  \ \ {\rm if } \ \ y\in (-1,0).
$$
This feedback control is optimal and it is
 consistent with the optimal open loop solution shown above.

 \begin{Remark}\label{Rem-opt-val-rep}
 {\rm If (\ref{e-opt-1}), (\ref{e-opt-2}) are valid, then the relationships (\ref{For-Remark}) are valid, the latter implying that
\begin{equation}\label{e-estimate}
{1\o T}(\bar\eta(y(T))-\bar\eta(y_0))=V(y_0)-{1\o T}\sum_{t=0}^{T-1} k(y(t),u(t)).
\end{equation}
This provides an interpretation of $\bar\eta(\cdot)$ as a function that defines the difference between the running cost $\ \disp {1\o T}\sum_{t=0}^{T-1} k(y(t),u(t))$ and  the optimal value $V(y_0)$ along the optimal trajectory.  Note that, if
$$
{1\o T}\sum_{t=0}^{T-1} k(y(t),u(t))= V_T(y_0),
$$
that is the process $(y(\cdot),u(\cdot)) $  is optimal on any finite time horizon as well, then (\ref{e-estimate}) can be rewritten as follows
\begin{equation}\label{e-estimate-2}
V_T(y_0) = V(y_0) - {1\o T} (\bar\eta(y(T))- \bar\eta(y_0))\ \ \ \forall \  T\geq 1.
\end{equation}
That was the case
 in the example considered above, in which the optimal trajectory $y(\cdot) $ satisfies the equalities: $\ y(T) = -y_0\ \forall y_0\in (0,1] $ and
 $\ y(T) = y_0\ \forall y_0\in [-1,0] $ for all $T\geq 1 $. This leads to $\ \bar\eta(y(T)) = 0 $ (see (\ref{e-example-11})) and, consequently, to that
 $$
 -{1\o T} (\bar\eta(y(T))- \bar\eta(y_0)) = {1\o T} \bar\eta(y_0) \ \ \ \forall \ T\geq 1 .
 $$
 Thus,
the relationships in (\ref{e-example-11-0}) are consistent with (\ref{e-estimate-2}). }
 \end{Remark}

{\bf Proof of Theorem \ref{ThN1}.} If  the pointwise limit (\ref{chesaro-lim-1-1}) exists, then, by Proposition \ref{Prop-imporant-1},
the limit function $V(\cdot) $  satisfies the inequality
$$
V(y_0)\geq d^*(y_0) \ \ \ \forall \ y_0\in Y.
$$
Therefore, to prove the statement (a), one needs to   show that
\begin{equation}\label{important-ineq-1}
 V(y_0)\leq d^*(y_0) \ \ \forall \ y_0\in Y.
\end{equation}
Similarly,  if  the pointwise limit (\ref{abel-lim-3-1}) exists, then, by Proposition \ref{Prop-imporant-1},
the limit function
 $h_{\a}(\cdot)$ satisfies  the inequality
$$
h(y_0)\geq d^*(y_0)\ \ \ \forall \ y_0\in Y.
$$
Therefore, to prove the statement (b), one needs to   show that
\begin{equation}\label{chesaro-lim-1-1-2-op}
h(y_0)\leq d^*(y_0)\ \ \ \forall \ y_0\in Y.
\end{equation}
{\it Proof of (\ref{important-ineq-1}).} Firstly, note that, by dividing (\ref{e-before-lim-1-1}) by $T$ and passing to the limit as $T\to \infty$, one obtains
\begin{equation}\label{V51-1}
V(y)\le V(f(y,u))  \ \ \forall \ (y,u)\in G.
\end{equation}
Also, by passing to the limit as $T\to \infty$ in (\ref{e-before-lim-1}), one obtains
\begin{equation}\label{V5-1}
\int_{G}V(y)\,\g(dy,du)\le \int_{G}k(y,u)\,\g(dy,du) \ \ \forall \ \g\in W.
\end{equation}
Inequality (\ref{V5-1})  can be rewritten in the form
\begin{equation*}
\int_{G}(k(y,u)-V(y))\,\gm\ge 0\quad \hbox{for all }\g\in W,
\end{equation*}
which is equivalent to that
\begin{equation}\label{CC15}
\min_{\g\in W}\int_{G}(k(y,u)-V(y))\,\gm\ge 0.
\end{equation}
The problem in the left hand side of the above inequality,
\begin{equation}\label{CC8}
\min_{\g\in W}\int_{G}(k(y,u)-V(y))\,\gm,
\end{equation}
is an IDLP problem, whose dual is
\begin{equation}\label{BB12}
\sup_{\eta\in C(Y)}\inf_{(y,u)\in G}\{k(y,u)-V(y)+\eta(f(y,u))-\eta(y)\}.
\end{equation}
Through equality of the optimal values of \eqref{CC8} and \eqref{BB12} (see Proposition 6 in \cite{GPS-2017}), we conclude that \eqref{CC15} is equivalent to
\begin{equation}\label{BB15}
\sup_{\eta\in C(Y)}\inf_{(y,u)\in G}\{k(y,u)-V(y)+\eta(f(y,u))-\eta(y)\}\ge 0.
\end{equation}
From \eqref{BB15} it follows that, for any $\ve>0$, there exists a function $\eta_{\ve}(\cdot)\in C(Y)$ such that
\begin{equation}\label{e-eps-feasib}
k(y,u)-V(y)+\eta_{\ve}(f(y,u))-\eta_{\ve}(y)\ge -\ve\quad \hbox{for all }(y,u)\in G.
\end{equation}
Consider the problem
\begin{equation}\label{BB21}
\sup_{(\psi,\eta)\in Q} \psi(y_0)= d^*(y_0),
\end{equation}
where
$Q$ is the set of pairs $(\psi,\eta)\in C(Y)\times C(Y)$ that satisfy inequalities
\begin{equation}\label{BB25}
\begin{aligned}
&k(y,u)-\psi(y)+\eta(f(y,u))-\eta(y)\ge 0,\\
&\psi(f(y,u))-\psi(y)\ge 0\quad \hbox{for all }(y,u)\in G.
\end{aligned}
\end{equation}
Note that the optimal value of problem (\ref{BB21}) is the same as that of (\ref{BB8}) (see (\ref{e-Pert-1}) in the proof of Lemma \ref{L-equivalince} taken with $\theta=0$).
Due to (\ref{V51-1}) and (\ref{e-eps-feasib}),  the pair $(\psi_{\ve}(\cdot), \eta_{\ve}(\cdot)) $, where  $\psi_{\ve}(\cdot):=V(\cdot)-\ve$, satisfies the inequalities \eqref{BB25}. Consequently,
$$
 d^*(y_0)\ge V(y_0)-\ve.
$$
This proves (\ref{important-ineq-1})
since $\ve > 0$ is arbitrarily small .

{\it Proof of (\ref{chesaro-lim-1-1-2-op}).} By passing to the limit  as $\a\up 1$  in (\ref{e-before-lim-2-1}),  we conclude that $h(\cdot)$ satisfies the inequality
 \begin{equation}\label{V51-2}
h(y)\le h(f(y,u)) \quad \hbox{for any }(y,u)\in G.
\end{equation}
Also, by  passing to the limit as $\a\up 1$ in (\ref{e-before-lim-2}) we establish that
\begin{equation}\label{V5-2}
\int_{G}h(y)\,\g(dy,du)\le \int_{G}k(y,u)\,\g(dy,du)  \ \ \forall \ \g\in W.
\end{equation}
Proceeding from this point  in exactly the same way as above, one establishes the validity of (\ref{chesaro-lim-1-1-2-op})
\hf


\section{Appendix }\label{LP-Duality-results}\label{Sec-Appendix}

\subsection{Another representation for the limit optimal values}\label{Sub-Sec-Another-Rep}

Let  $\K$ be the set of continuous functions that satisfy the following relationships:
\begin{equation}\label{V51}
w(y)\le w(f(y,u)) \quad \hbox{for any }(y,u)\in G
\end{equation}
and
\begin{equation}\label{V5}
\int_{G}w(y)\,\g(dy,du)\le \int_{G}k(y,u)\,\g(dy,du)\quad \hbox{for all }\g\in W.
\end{equation}
In these notations,  the relationships (\ref{V51-1}), (\ref{V5-1}) and  (\ref{V51-2}), (\ref{V5-2})  are equivalent to the inclusions
\begin{equation}\label{e-Marc}
V(\cdot)\in \K,
\end{equation}
and
\begin{equation}\label{e-Marc-4}
 h(\cdot)\in \K,
\end{equation}
respectively.

\begin{Proposition}\label{Th-C}
 {\rm (a)} Let the pointwise limit (\ref{chesaro-lim-1-1}) exist and the function
 $V (\cdot)  $ be continuous.  Then
\begin{equation}\label{chesaro-lim-2}
 V(y_0)= \sup\{w(y_0)\ | \ w(\cdot)\in \K\} \ \ \forall \ y_0\in Y.
\end{equation}
{\rm (b)} Let the pointwise limit (\ref{abel-lim-3-1})  exists and
the function $h (\cdot)  $ be continuous. Then
\begin{equation}\label{abel-lim-4}
 h(y_0)= \sup\{w(y_0)\ | \ w(\cdot)\in \K\} \ \ \forall \ y_0\in Y.
\end{equation}
\end{Proposition}
{\bf Proof.} Note that, due to (\ref{e-Marc}) and (\ref{e-Marc-4})
\begin{equation}\label{e-Marc-1}
V(y_0)\leq \sup\{w(y_0)\ | \ w(\cdot)\in \K\},  \ \ \ \ \ h(y_0)\leq \sup\{w(y_0)\ | \ w(\cdot)\in \K\}\ \ \ \forall \ y_0\in Y.
\end{equation}
Therefore, to prove the proposition, it is sufficient to establish that the inequalities opposite  to (\ref{e-Marc-1})  are valid.
For a natural $T$, let $u_T(\cdot)$ be an optimal  control in \eqref{A112-1}, $\g_T\in \G_T(y_0)$ be the occupational measure generated by this control, and $y_T(\cdot)$ be the corresponding trajectory. Then
$$
{1\o T}\sum_{t=0}^{T-1} k(y_T(t),u_T(t))=\int_G k(y,u)\,\g_T(dy,du)=V_T(y_0).
$$
Let $\g_T(dy,du)$ converge to $\g$ in weak$^*$ topology as $T\to \infty$ along a subsequence (we do not relabel). Note that $\g\in W $ (due to (\ref{convergence-to-W})). From the equality above, by passing to the limit as $T\to \infty$, we obtain
\begin{equation}\label{V6}
\int_{G}k(y,u)\,\g(dy,du)=V(y_0).
\end{equation}
For $w\in \K$, taking into account the monotonicity property \eqref{V51}, we have
$$
w(y_0)={1\o T}\sum_{t=0}^{T-1} w(y_0) \le {1\o T}\sum_{t=0}^{T-1} w(y_T(t))=\int_{G} w(y)\,\g_T(dy,du).
$$
Since $w$ is continuous, we can pass to the limit as $T\to\infty$ and obtain
$$
w(y_0)\le \int_{G} w(y)\,\g(dy,du).
$$
Combining this with \eqref{V5} and \eqref{V6} we obtain
$$
w(y_0)\le \int_{G} w(y)\,\g(dy,du)\le \int_{G} k(y,u)\,\g(dy,du)= V(y_0).
$$
The latter implies that the inequality opposite to the first inequality in (\ref{e-Marc-1}) is valid. This  proves part (a) of the proposition.

The proof of the inequality opposite to  the second inequality in (\ref{e-Marc-1})  is similar. For $\a\in (0,1), $ let $u_{\a}(\cdot)$ be an optimal  control in \eqref{A112-2}, $\g_{\a}\in \T_{\a}(y_0)$ be the occupational measure generated by this control, and $y_{\a}(\cdot)$ be the corresponding trajectory. Then
$$
(1-\a)\sum_0^{\infty}\a^t k(y_{\a}(t),u_{\a}(t))=\int_G k(y,u)\,\g_{\a}(dy,du)=h_{\a}(y_0).
$$
Let $\g_{\a}(dy,du)$ converge to $\g$ in weak$^*$ topology as $\a\to 1$ along a subsequence (we do not relabel). Note that $\g\in W $ (due to (\ref{convergence-to-W})). From the equality above, by passing to the limit as $\a\to 1$ we obtain
\begin{equation}\label{V6-a}
\int_{G}k(y,u)\,\g(dy,du)=h(y_0).
\end{equation}
Combining this with \eqref{V5} and \eqref{V6-a} we obtain
$$
w(y_0)\le \int_{G} w(y)\,\g(dy,du)\le \int_{G} k(y,u)\,\g(dy,du)= h(y_0).
$$
The latter implies that the inequality opposite to the second inequality in (\ref{e-Marc-1}) is valid, and, thus,  proves part (b) of the proposition.
\hf
\begin{Remark}\label{Rem-LP-Rep}
 {\rm It can be verified directly that the optimal value of the problem in the  right hand side of (\ref{chesaro-lim-2}) and (\ref{abel-lim-4}) is equal to $d^*(y_0) $ (the optimal value of the dual problem (\ref{BB8})). Results establishing the validity of presentations similar to (\ref{chesaro-lim-2}) and (\ref{abel-lim-4}) in continuous time setting were obtained in  \cite{BQR-2015}. }
\end{Remark}

\subsection{Results referred to in Sections \ref{Sec-Up-Bound} and \ref{Sec-equality}}\label{Sub-Sec-Auxiliary-Res}
Consider a perturbed version of the IDLP problem (\ref{BB1})
\begin{equation}\label{BB1-Per}
\inf_{(\g,\xi)\in \O(y_0)} \left\{\int_G k(y,u)\gm + \theta \int_G \xi(dy,du)\right\}  := k^*(\theta, y_0),
\end{equation}
and the corresponding perturbed version of the dual problem
(\ref{BB8})
\begin{equation}\label{BB8-Per}
\sup_{(\mu,\psi,\eta)\in \D(\theta,y_0)} \mu=:d^*(\theta, y_0),
\end{equation}
where $\D(\theta,y_0)$ is the set of triplets $(\mu,\psi(\cdot),\eta(\cdot))\in \reals\times C(Y)\times C(Y)$  that satisfy the inequalities
\begin{equation}\label{BB7-Per}
\begin{aligned}
&k(y,u)+(\psi(y_0)-\psi(y))+\eta(f(y,u))-\eta(y)-\mu \ge 0,\\
&\psi(f(y,u))-\psi(y)\ge -\theta \ \ \ \ \forall \ (y,u)\in G,
\end{aligned}
\end{equation}
Note that $\theta \geq 0 $ is a perturbation parameter and note that (\ref{BB1-Per}) and (\ref{BB8-Per}) become (\ref{BB1}) and (\ref{BB8}) with $\theta = 0 $.
 Consider also the problem
\begin{equation}\label{BB21-1-Per}
\sup_{(\psi,\eta)\in Q(\theta)} \psi(y_0):= \bar d^*(\theta,y_0),
\end{equation}
where
$Q(\theta)$ is the set of pairs $(\psi(\cdot),\eta(\cdot))\in C(Y)\times C(Y)$ that satisfy the inequalities
\begin{equation}\label{BB25-1-Per}
\begin{aligned}
&k(y,u)-\psi(y)+\eta(f(y,u))-\eta(y)\ge 0,\\
&\psi(f(y,u))-\psi(y)\ge -\theta\ \ \ \ \forall \ (y,u)\in G,
\end{aligned}
\end{equation}

\begin{Lemma}\label{L-equivalince}
The following relationships are valid:
\begin{equation}\label{e-Pert-1}
\bar d^*(\theta,y_0) = d^*(\theta,y_0)\leq k^*(\theta, y_0)  \ \ \ \forall \ \theta\geq 0.
\end{equation}
\end{Lemma}
{\bf Proof.} Let us prove, first, that
\begin{equation}\label{e-Pert-2}
\bar d^*(\theta,y_0) = d^*(\theta,y_0)  \ \ \ \forall \ \theta\geq 0.
\end{equation}
In fact, the inequality $\bar d^*(\theta,y_0) \leq  d^*(\theta,y_0) $ is true (since, for any pair $\ (\psi(\cdot),\eta(\cdot))\in Q(\theta) $, the triplet $\ (\mu , \psi(\cdot),\eta(\cdot))\in \D(\theta,y_0) $ with $\mu = \psi(y_0) $). Let us prove the opposite inequality. Let a triplet $\ (\mu' , \psi'(\cdot),\eta'(\cdot))\in \D(\theta,y_0) $ be such that $\mu'\geq d^*(\theta,y_0)-\delta $, with $\delta > 0 $ being arbitrarily small. Then the pair $\ (\tilde \psi'(\cdot), \eta'(\cdot))\in Q(\theta) $, with $\tilde\psi'(y)= \psi'(y)- \psi'(y_0) + \mu'$. Since $ \tilde\psi'(y_0) = \mu'$, it leads to the inequality $\bar d^*(\theta,y_0) \geq  d^*(\theta,y_0) - \delta $ and, consequently, to the inequality $\bar d^*(\theta,y_0) \geq  d^*(\theta,y_0) $ since $\delta > 0 $ is arbitrarily small. Thus, (\ref{e-Pert-2}) is proved.

Let us now prove the inequality
\begin{equation}\label{e-Pert-3}
d^*(\theta,y_0)\leq k^*(\theta, y_0)  \ \ \ \forall \ \theta\geq 0.
\end{equation}
Take any $(\g,\xi)\in \O(y_0)$ and $(\mu,\psi,\eta)\in \D(\theta, y_0)$. Integrating the first inequality in \eqref{BB7-Per} with respect to $\g$ and taking into account that $\g\in W$ we conclude that
$$
\int_G k(y,u)\gm+\int_G(\psi(y_0)-\psi(y))\gm\ge \mu.
$$
Taking into account that $(\g,\xi)\in \O(y_0)$ and the second inequality in \eqref{BB7-Per}, we obtain
$$
\int_G(\psi(y_0)-\psi(y))\gamma(dy,du)=-\int_G (\psi(f(y,u))-\psi(y))\xi(dy,du)\le \theta\int_{G}\xi(dy,du).
$$
Therefore,
$$
\int_G k(y,u)\gm + \theta\int_{G}\xi(dy,du) \ge \mu .
$$
This proves (\ref{e-Pert-3}). \hf

Let $C^*(Y) $ stand for the space of continuous linear functionals  on $C(Y) $ and let $\mathcal{M}(G) $ stand for the space of measures defined on Borel subsets of $G$. Define a linear operator $\mathcal{A}(\cdot): \mathcal{M}(G) \times  \mathcal{M}(G)\mapsto\reals^1\times  C^*(Y) \times  C^*(Y) $  as follows: for any $(\gamma, \xi)\in \M(G)\times \M(G) $,
\begin{equation}\label{e-Duality-1}
\mathcal{A}(\gamma, \xi):= \left(\int_{G}\gamma(dy,du), \ a_{(\gamma, \xi)}, \ b_{\gamma}\right),
\end{equation}
where $a_{(\gamma, \xi)}, \ b_{\gamma}\in C^*(Y)$ are defined by the equation: $\ \forall \  \phi(\cdot)\in C(Y)$,
$$
 \ a_{(\gamma, \xi)}(\phi) := -\left\{\int_{G} (\phi(y_0)-\phi(y))\gamma(dy,du)
 + \int_{G}(\phi(f(y,u)) - \phi(y))\xi(dy,du)\right\},
 $$
 \vspace{-0.4cm}
$$
\  b_{\gamma}(\phi):= - \left\{\int_{G}(\phi(f(y,u)) - \phi(y))\gamma(dy,du)\right\}.
$$
In this notation, the set $\Omega(y_0)$ defined in (\ref{eq-Omega}) can be rewritten as follows
$$
 \Omega(y_0)= \{(\gamma , \xi)\in \M_+(G)\times \M_+(G)\ : \ \A(\gamma, \xi)= (1, {\bf 0}, {\bf 0})\},
 $$
where  ${\bf 0} $ stands for the zero element of $C^*(Y)$. Also,
problem (\ref{BB1}) takes the form
\begin{equation}\label{e-Duality-2-0}
\inf_{(\gamma , \xi)\in \Omega(y_0)}\langle k, \gamma  \rangle\ = k^*(y_0),
\end{equation}
where
 $\langle \cdot, \gamma  \rangle $ (also,
$\langle \cdot, \xi  \rangle $ in the sequel) denoting the integral of the corresponding function over $\gamma$ (respectively, over $\xi$).
Note that, for any $(\mu, \psi(\cdot), \eta(\cdot))\in \reals^1\times C(Y)\times C(Y) $,
$$
\langle A(\gamma, \xi), (\mu, \psi, \eta)\rangle =
\mu \int_{G}\gamma(dy,du) + a_{(\gamma, \xi)}(\psi)
+ b_{\gamma}(\eta)
$$
$$
=  \int_{G}\left(\mu - (\psi(y_0)-\psi(y)) - ( \eta(f(y,u)) - \eta(y))\right)\gamma(dy,du)
$$
$$
- \int_{G} ( \psi(f(y,u)) - \psi(y))\xi(dy,du).
$$
Define now the linear operator\\ $\A^*(\cdot): \reals^1\times C(Y)\times C(Y) \mapsto C(G)\times C(G)\subset \M^*(G)\times  \M^*(G) $ in such a way that, for any $(\mu, \psi(\cdot),\eta(\cdot))\in \reals^1\times C(Y)\times C(Y) $,
$$
\A^*(\mu, \psi, \eta)(y,u):= \left(\mu - (\psi(y_0)-\psi(y)) - ( \eta(f(y,u)) - \eta(y)), \ -( \psi(f(y,u)) - \psi(y))\right).
$$
Thus,
$$
\langle  \A^*(\mu, \psi, \eta), (\gamma, \xi) \rangle =  \int_{G}\left(\mu - (\psi(y_0)-\psi(y)) - ( \eta(f(y,u)) - \eta(y))\right)\gamma(dy,du)
$$
$$
 - \int_{G} ( \psi(f(y,u)) - \psi(y))\xi(dy,du) = \langle A(\gamma, \xi), (\mu, \psi, \eta)\rangle  .
$$
That is, the operator $ \A^*(\cdot) $ is the adjoint of  $ \A(\cdot) $. The problem dual to (\ref{e-Duality-2-0})
is of the form (see \cite{And-1} and \cite{And-2})
$$
\sup_{(\mu, \psi(\cdot), \eta(\cdot))\in \reals^1\times C(Y)\times C(Y)} \mu = d^*(y_0)
$$
$$
\ s.\ t.  \ \ \ \ \ \ \ \ \ \ \   \ \ \ \ \ \ \ \ \ \ \ \ \ \ \ \ \ \ \ \   \ \ \ \ \ \ \ \ \
$$
$$
 - \A^*(\mu, \psi, \eta)(y,u) + (k(y,u), 0)\geq (0,0) \ \ \forall (y,u)\in G,
 $$
the latter being equivalent to (\ref{BB8}).

{\bf Proof of Lemma \ref{L2}.}
Let
$$
H:= \Big\{\left( \A(\gamma, \xi), \int_{G}k(y,u)\gamma(dy,du) + r\right)\ :
$$
\vspace{-0.4cm}
$$
\ (\gamma, \xi)\in \M_+(G)\times \M_+(G), \ r\geq 0\Big\}\subset \reals^1\times C^*(Y)\times C^*(Y)\times \reals^1  ,
$$
and let $\bar{H}$ stand for the closure of $H$ in the weak$^*$ topology of $\reals^1\times C^*(Y)\times C^*(Y) \times \reals^1$. Consider the problem
\begin{equation}\label{limits-non-ergodic-pert-dual-4-1}
\inf \{\theta \ | \ (1, {\bf 0}, {\bf 0}, \theta)\in \bar{H} \}:= k_{sub}^*(y_0).
\end{equation}
Its optimal value $k_{sub}^*(y_0)$ is called the subvalue of the IDLP problem (\ref{e-Duality-2-0}). Let us show that the optimal value of (\ref{CC14}) is equal to the subvalue. In fact,
as can be readily seen, $\left(1, {\bf 0}, {\bf 0},  \int_{G}k(y,u)\gamma(dy,du)\right)\in \bar{H} $  if $\gamma\in W_2(y_0)$. Consequently,
$$
k_{sub}^*(y_0)\leq \min_{\gamma\in W_2(y_0)}\int_{G}k(y,u)\gamma(dy,du).
$$
From the fact that $k_{sub}^*(y_0)$ is defined as  the optimal value in (\ref{limits-non-ergodic-pert-dual-4-1}) it follows that there exists  a sequence $(\gamma_l,\xi_l)\in \M_+(G)\times \M_+(G) $  such that $\A(\gamma_l,\xi_l) $ converges (in weak$^*$ topology)
to $(1, {\bf 0}, {\bf 0})$, with   $\int_{G}k(y,u)\gamma_l(dy,du) $ converging to $k_{sub}^*(y_0)$ as $l$ tends to infinity. That is (see (\ref{e-Duality-1})),
$$
\int_{G}\gamma_l(dy,du)\rightarrow 1, \ \ a_{(\gamma_l, \xi_l)}\rightarrow {\bf 0},\ \  b_{\gamma_l}\rightarrow {\bf 0},
$$
 \vspace{-0.4cm}
$$
\ \ \int_{G}k(y,u)\gamma_l(dy,du)\rightarrow k_{sub}^*(y_0).
$$
Without loss of generality, one may assume
that $\gamma_l$ converges in weak$^*$ topology
 to a  measure $\gamma$ that satisfies the relationships
 $$
 \int_{G}\gamma(dy,du)=1,  \ \  b_{\gamma}= {\bf 0}\ \ \ \ \Rightarrow \ \ \ \  \gamma\in W.
 $$
 Also,
 $
 \  a_{(\gamma, \xi_l)}\rightarrow {\bf 0}
 $ and $\int_{G}k(y,u)\gamma(dy,du)= k_{sub}^*(y_0) $. That is, $\gamma\in W_2(y_0)$ and therefore,
 $$
\min_{\gamma\in W_2(y_0)}\int_{G}k(y,u)\gamma(dy,du)\leq k_{sub}^*(y_0).
$$
 Thus, the optimal value of (\ref{CC14}) is equal to the subvalue.
 To complete the proof, it is sufficient to note that the subvalue of an IDLP problem is equal to the optimal value of its dual provided that the former is bounded (see, e.g., Theorem 3 in \cite{And-1}). That is, $k_{sub}^*(y_0) = d^*(y_0)$. \hf

Let us conclude this section with proving the validity of the following proposition.

\begin{Proposition}\label{PN3} The optimal value of the problem  in the left hand side of (\ref{CC21}) is equal to $\liminf_{T\to \infty} V_T(y_0) $. That is,
$$
\inf_{u(\cdot)\in \U(y_0)} \liminf_{T\to \infty} {1\o T}\sum_{t=0}^{T-1} k(y(t),u(t))=\liminf_{T\to \infty} V_T(y_0)\ \ \ \ \forall \ y_0\in Y.
$$
\end{Proposition}

{\bf Proof.} 
Let $u(\cdot)\in \U(y_0)$ and let $y(\cdot)$ be the corresponding trajectory. Then
$$
{1\o T}\sum_{t=0}^{T-1} k(y(t),u(t))\ge V_T(y_0).
$$
Therefore,
$$
\liminf_{T\to \infty}{1\o T}\sum_{t=0}^{T-1} k(y(t),u(t))\ge \liminf_{T\to \infty} V_T(y_0)
$$
and, hence,
$$
\inf_{u(\cdot)\in \U(y_0)} \liminf_{T\to \infty} {1\o T}\sum_{t=0}^{T-1} k(y(t),u(t))\ge\liminf_{T\to \infty} V_T(y_0).
$$
Let us prove the opposite inequality.
For any $\varepsilon>0$ and $u(\cdot)\in \U(y_0)$,  and for sufficiently large $T$,
$$
{1\o T}\sum_{t=0}^{T-1} k(u(t),y(t)) \ge \liminf_{T\to \infty} {1\o T}\sum_{t=0}^{T-1} k(u(t),y(t)) -\varepsilon ,
$$
where $y(\cdot)=y(t,y_0,u) $.
Therefore,
\begin{equation*}
\begin{aligned}
{1\o T}\sum_{t=0}^{T-1} k(u(t),y(t)) \ge \inf_{u'\in \U(y_0)}\liminf_{T\to \infty} {1\o T}\sum_{t=0}^{T-1} k(u'(t),y'(t)) -\varepsilon
\end{aligned}
\end{equation*}
and, consequently,
\begin{equation*}
\begin{aligned}
V_T(y_0) \ge \inf_{u\in \U(y_0)}\liminf_{T\to \infty} {1\o T}\sum_{t=0}^{T-1} k(u(t),y(t)) -\varepsilon.
\end{aligned}
\end{equation*}
Hence,
$$
\liminf_{T\to \infty}V_T(y_0)\ge \inf_{u\in \U(y_0)}\liminf_{T\to \infty} {1\o T}\sum_{t=0}^{T-1} k(u(t),y(t)).
$$
The proposition is proved. \hf

\section{Conclusions}\label{Sec-COnclusions}

\ \ \ \ We have introduced the IDLP problem, the optimal value of which gives an upper bound for  $\ \limsup_{T\rightarrow\infty}V_T(y_0)$ and $\ \limsup_{\a\up 1}h_{\a}(y_0)$, with the optimal value of the corresponding dual problem providing
 a lower bound for $\ \liminf_{T\rightarrow\infty}V_T(y_0)$ and $\ \liminf_{\a\up 1}h_{\a}(y_0)$. While the result establishing the validity of the lower bound (Proposition \ref{Prop-imporant-1}) is very similar to the corresponding result in \cite{BG-2018}, the statement about the validity of the upper bound (Theorem \ref{Th-upper-bound}) is much stronger than its continuous time counterpart in \cite{BG-2018}, where it was assumed that
  the uniform limits $\ \lim_{T\rightarrow\infty}V_T(y_0)$ and $\ \lim_{\a\up 1}h_{\a}(y_0)$  exist and are Lipschitz continuous. Note also that,  in contrast to the result of \cite{BG-2018},  we did not assume that the set $Y$ is  invariant (only that it is viable).   We believe that establishing the validity of the upper bound  for systems evolving in continuous time under assumptions similar to those of  Theorem \ref{Th-upper-bound} is possible, and it   can be a subject for future research.

  We have also established that, if the pointwise limits $\ \lim_{T\rightarrow\infty}V_T(y_0)$ and $\ \lim_{\a\up 1}h_{\a}(y_0)$ exist and are continuous, then they are equal to the optimal value of the dual problem (Theorem \ref{ThN1}). A similar statement in the continuous time setting can be established  using a similar argument if the limits of the optimal value functions exist and are continuously differentiable. This assumption is, however, too strong, and finding less restrictive conditions, under which a statement similar to Theorem \ref{ThN1} for systems in continuous time is valid, can also be a subject for future research.

Finally, we have stated sufficient and necessary optimality conditions for the long-run average optimal control problem  using an optimal solution of the dual problem (Propositions \ref{Prop-optim-cond-suf} and \ref{Prop-optim-cond-nec}). Similar results can be  readily  obtained in the continuous time case too.

\bigskip

{\bf Acknowledgment.} We would like  to express our gratitude to D. Khlopin and to M. Quincampoix
 for useful discussions and for sharing with us some insightful examples.

\end{document}